\documentclass[12pt]{amsart}

\DeclareMathOperator{\Aut}{Aut}%
\DeclareMathOperator{\chr}{char}%
\DeclareMathOperator{\Gal}{Gal}%
\DeclareMathOperator{\pr}{pr}%
\DeclareMathOperator{\rank}{rank}

\usepackage{amssymb}
\usepackage{hyperref}

\begin{document}

\title[Galois Module Structure of $p$th-Power Classes]{Galois
Module Structure of $p$th-Power Classes of Cyclic Extensions of
Degree $p^n$}

\dedicatory{Dedicated to the memory of Walter Feit}

\author[J\'{a}n Min\'{a}\v{c}]{J\'an Min\'a\v{c}}
\address{Department of Mathematics, Middlesex College, \ University
of Western Ontario, London, Ontario \ N6A 5B7 \ CANADA}
\email{minac@uwo.ca}

\author[Andrew Schultz]{Andrew Schultz}
\address{Department of Mathematics, Building 380, Stanford
University, Stanford, California \ 94305-2125 \ USA}
\email{aschultz@stanford.edu}

\author[John Swallow]{John Swallow}
\address{Department of Mathematics, Davidson College, Box 7046,
Davidson, North Carolina \ 28035-7046 \ USA}
\email{joswallow@davidson.edu}

\begin{abstract}
In the mid-1960s Borevi\v{c} and Faddeev initiated the study of the
Galois module structure of groups of $p$th-power classes of cyclic
extensions $K/F$ of $p$th-power degree. They determined the
structure of these modules in the case when $F$ is a local field. In
this paper we determine these Galois modules for all base fields
$F$.
\end{abstract}



\thanks{J\'an Min\'a\v{c} was supported in part by National Sciences
and Engineering Research Council of Canada grant R0370A01, by a
Distinguished Professorship for 2004--2005 at the University of
Western Ontario, and by the Mathematical Sciences Research
Institute, Berkeley.  John Swallow was supported in part by National
Security Agency grant MDA904-02-1-0061.}

\maketitle

\newtheorem{theorem}{Theorem}
\newtheorem{corollary}{Corollary}
\newtheorem{proposition}{Proposition}
\newtheorem{lemma}{Lemma}

\newtheorem*{definition*}{Definition}
\newtheorem*{remark*}{Remark}

\newcommand{\C}{\mathbb{C}}
\newcommand{\Ec}{\mathcal{E}}
\newcommand{\Fp}{\mathbb{F}_p}
\newcommand{\Ft}{\mathbb{F}_2}
\newcommand{\Gc}{\mathcal{G}}
\newcommand{\Ic}{\mathcal{I}}
\newcommand{\N}{\mathbb{N}}
\newcommand{\R}{\mathbb{R}}
\newcommand{\Z}{\mathbb{Z}}
\newcommand{\Zc}{\mathcal{Z}}

In 1947 \v{S}afarevi\v{c} initiated the study of Galois groups of
maximal $p$-extensions of fields with the case of local
fields~\cite{Sa}, and this study has grown into what is both an
elegant theory as well as an efficient tool in the arithmetic of
fields. From the very beginning it became clear that the groups of
$p$th-power classes of the various field extensions of a base field
encode basic information about the structure of the Galois groups of
maximal $p$-extensions. (See \cite{Koc} and \cite{Ser}.) Such groups
of $p$th-power classes arise naturally in studies in arithmetic
algebraic geometry, for example in the study of elliptic curves.

In 1960 Faddeev began to study the Galois module structure of
$p$th-power classes of cyclic $p$-extensions, again in the case of
local fields, and during the mid-1960s he and Borevi\v{c}
established the structure of these Galois modules using basic
arithmetic invariants attached to Galois extensions. (See \cite{F}
and \cite{Bo}.) In 2003 two of the authors ascertained the Galois
module structure of $p$th-power classes in the case of cyclic
extensions of degree $p$ over all base fields $F$ containing a
primitive $p$th root of unity \cite{MS}. Very recently, this work
paved the way for the determination of the entire Galois cohomology
as a Galois module in the case of a cyclic extension of degree $p$
of a base field containing a primitive $p$th root of unity, using
Voevodsky's recent work on Galois cohomology (\cite{LMS}; see
\cite{V1} and \cite{V2}).

In this paper we extend the results obtained in \cite{MS} in two
directions. First, our results hold for cyclic extensions of any
$p$th-power degree, rather than just $p$, and, furthermore, we no
longer require that the base field contain a primitive $p$th root
of unity. Thus our results provide a complete determination of
$p$th-power classes as Galois modules for all cyclic extensions of
$p$th-power degree.

We expect that, just as the results and techniques in \cite{MS}
helped to determine the entire Milnor $K$-theory modulo $p$ as a
Galois module in the case of cyclic extensions of degree $p$, so
will the results and methods developed in this paper lead to the
determination of the entire Milnor $K$-theory modulo $p$ as a
Galois module in the case of cyclic extensions of $p$th-power
degree. In fact, precisely such a generalization has already taken
place in the case of characteristic $p$ \cite{BLMS}.

Similarly, in the same way as the results and techniques developed
in \cite{MS} led in \cite{MS2} to the solution of Galois embedding
problems and the discovery of a new automatic realization of Galois
groups, it is clear that the results in this paper will also have
such Galois-theoretic applications. In a subsequent paper we will
consider some of these applications.

Our basic approach to the problem is induction, and some of the
results in \cite{MS} handle the base case.  In the end, however,
neither the results nor the techniques employed are
straightforward generalizations of the work in \cite{MS}. First,
the possible generalization of the innocent summand of dimension
$1$ or $2$ considered in \cite{MS} turned out to be rather subtle to
handle. These new summands of dimension $p^i+1$ for some
$i\in\N$ are very interesting invariants of cyclic extensions of
$p$th-power degree. Another substantial challenge was to generate
enough norms, and the resolution involves several thorny induction
arguments. Finally, the case $p=2$ presented a new problem for
quartic extensions, and this problem is taken care of as a
separate base induction case.

Fundamentally, the classification of $p$th-power classes as Galois
modules depends upon arithmetic invariants, all of which originate
from the images of the norms of the intermediate fields of $K/F$.
The classification, in short, has the flavor of local class field
theory, and although the arguments underlying the classification
are not straightforward, the final results, just as in local class
field theory, have a rather simple and elegant form, which we now
describe.

Let $p$ be a prime number, $n\ge 1$ an integer, $F$ an arbitrary
field, and $K$ a Galois extension of $F$ with group $G = \langle
\sigma \rangle$ cyclic of order $p^n$. Let $F^\times$ denote the
multiplicative group of nonzero elements of $F$.  Let $J = J(K) =
K^\times/K^{\times p}$ be the $\Fp[G]$-module of $p$th-power
classes, where for a given $\gamma\in K^\times$ we write the class
$\gamma K^{\times p}$ as $[\gamma]$.  Similarly, let $J(F) =
F^\times/F^{\times p}$ be the $\Fp$-module of $p$th-power classes of
$F^\times$, where for a given $f\in F^\times$ we write the class $f
F^{\times p}$ as $[f]_F$.  Let $N_{K/F}\colon K\to F$ be the norm
map, and write $N\colon K^\times/K^{\times p}\to F^\times/F^{\times
p}$ for the map induced by $N_{K/F}$. Also by abuse of notation we
use the same symbol $N$ to denote the endomorphism $N\colon K^\times
/K^{\times p} \to K^\times/K^{\times p}$ induced by $N\colon
K^\times/K^{\times p} \to F^\times/F^{\times p}$ defined above,
followed by the map induced by the inclusion map $\epsilon =
\epsilon_K \colon F^\times \to K^ \times$.

Further let $K_i$, $i=0, \dots, n$, be the intermediate field of
$K/F$ such that $[K_i:F]=p^i$. Denote by $H_i$ the Galois group
$\Gal(K/K_i)\subset G$.  Let $[K_i^\times]$ denote the submodule of
$J$ which is the image of the map induced by the inclusion map
$K_i^\times\to K^\times: \left[K_i^\times\right] = K_i^\times
K^{\times p}/K^{\times p}$.  Similarly, for other $G$-submodules
$A\subset K^\times$, such as $A=N_{K_i/F}(K_i^\times)$, let $[A] =
AK^{\times}/K^{\times p}$.

If the characteristic of $F$ is not $p$, we denote by $\xi_p$ a
fixed primitive $p$th root of unity in a fixed algebraic closure of
$F$.  We say that the condition $\xi_p\not\in F$ is satisfied if
either the characteristic of $F$ is $p$ or the characteristic of $F$
is not $p$ and $\xi_p\not\in F$.

\begin{theorem}\label{th:nox}
    Suppose that either
    \begin{itemize}
        \item $\xi_p\not\in F$, or
        \item $p=2$, $n=1$, and $-1\not\in N_{K/F}(K^\times)$.
    \end{itemize}
    Then the $\Fp[G]$-module $J$ decomposes as
    \begin{equation*}
        J = Y_n \oplus Y_{n-1} \oplus \dots \oplus Y_0,
    \end{equation*}
    where $Y_i$ is a direct sum of cyclic $\Fp[G]$-modules
    of dimension $p^i$ and
    \begin{equation*}
        [K_i^\times] = J^{H_i}, \quad 0\le i\le n.
    \end{equation*}
\end{theorem}

It is easy to show that this decomposition of $J$ is unique up to
isomorphism. (In fact this also follows from a well-known result of
Azumaya. See \cite[page~144]{AnFu}.) In the following corollary we
determine the sizes of the modules $Y_i$ in terms of norms. Observe
that direct sums of cyclic $\Fp[G]$-modules of dimension $p^i$ are
free $\Fp[G/H_ i]$-modules.  Let
\begin{equation*}
        e_i = \dim_{\Fp}
        \left(\left[N_{K_i/F}\left(K_i^\times\right)
    \right]/\left[N_{K_{i+1}/F}\left(K_{i+1}^\times
    \right)\right]\right), \quad 0\le i<n,
\end{equation*}
and let $e_n= \dim_{\Fp} [N_{K/F}(K^\times)]$.

\begin{corollary}\label{co:nox}
    For each $0\le i\le n$
    \begin{equation*}
        [N_{K_i/F}(K_i^\times)] = (Y_i+Y_{i+1}+\dots+Y_n)^G,
    \end{equation*}
    and
    \begin{equation*}
        \rank_{\Fp[G/H_i]} Y_i = e_i.
    \end{equation*}
\end{corollary}

For $K/F$ not satisfying the conditions of the theorem above, we
adopt the conventions $K_{-\infty}^\times = \{1\}$ and
$p^{-\infty}=0$ and make the following definition.

\begin{definition*}[Exceptional Element]
    Suppose that $\xi_p\in F$ and, if $p=2$, that either $n>1$ or
    $-1\in N_{K/F}(K^\times)$.  We set
    \begin{align*}
        i(K/F) := \min \{\ &i\in \{-\infty, 0, 1, \dots, n\} \
        \vert \ \exists \delta\in K^\times \text{ such that } \\
        &[N_{K/F}(\delta)]_F \neq [1]_F \text{\ and\ } \\
        &[\delta]^{\tau-1} \in [K_i^\times]\ \forall
        \tau\in\Gal(K/F)\}.
    \end{align*}
    We say that $\delta\in K^\times$ is an \emph{exceptional
    element} of $K/F$ if $[N_{K/F}(\delta)]_F \neq [1]_F$ and
    $[\delta]^{\tau-1}\in [K_{i(K/F)}^\times]$ for all $\tau
    \in \Gal(K/F)$. Elements of $K^\times$ that are not
    exceptional are said to be \emph{unexceptional}.  For
    simplicity, we often write $m$ instead of $i(K/F)$.
\end{definition*}

Observe that $[\delta]^{(\tau-1)}\in [K_i^\times]$ for all
$\tau\in G$ if and only if $[\delta]^{(\sigma-1)}\in
[K_i^\times]$ for a fixed generator $\sigma\in G$.  In what
follows we will use this formulation for our given generator
$\sigma$.

Note that if $\delta$ is an exceptional element then
$m=i(K/F)=-\infty$ if and only if $[\delta]^\sigma = [\delta]$ and
$[N_{K/F}(\delta)]_F \neq [1]_F$.

Because the exceptionality of an element $\gamma\in K^\times$ is
independent of the particular representative $\gamma$ of
$[\gamma]$, we define $[\gamma]$ to be exceptional if $\gamma$ is
exceptional.  It is also useful to observe that if an
$\Fp[G]$-generator $[\gamma]$ of a module $M_\gamma\subset J$ is
exceptional, then so is any other $\Fp[G]$-generator $[\omega]$ of
$M_\gamma$.  Indeed, using additive notation for $J$ for the
moment, any such generator $[\omega]$ has the form
\begin{equation*}
    [\omega] = c_0[\gamma] + c_1(\sigma-1)[\gamma] +
    c_2(\sigma-1)^2[\gamma] + \dots, \quad c_0, c_1, \dots\in \Fp,
    \quad c_0\neq 0.
\end{equation*}
Then $[N_{K/F}(\omega)]_F = [N_{K/F}(\gamma)]^{c_0}_F \neq [1]$
and $[\omega]^{\sigma-1} \in [K_{m}^\times]$.

In Proposition~\ref{pr:excexist} we show that exceptional elements
always exist for $K/F$ satisfying the hypothesis in the Definition
above, and in Proposition~\ref{pr:ldelta} we show that, in fact,
$m\le n-1$. Finally, note that since $N_{K/F}(K_{n-1}^\times)
\subset F^{\times p}$, each exceptional element $\delta \in
K_n^\times \setminus K_{n-1}^\times$.

Moreover, for these $K/F$, we have Kummer theory, because
$\xi_p\in F$. Hence $K_1=F(\root{p}\of{a})$ for some $a\in F$.  In
section~\ref{se:exc} we prove some more specific results about
exceptional elements in terms of $a$: exceptional elements satisfy
$[N_{K/F}(\delta)]_F = [a]_F^s$ for $s\not\equiv 0\bmod p$ and
that for all $K/F$ as above, an exceptional element $\delta\in
K^\times$ exists satisfying $[N_{K/F}(\delta)]_F = [a]_F$.

\begin{theorem}\label{th:x}
    Suppose that $\xi_p\in F$ and, if $p=2$, that either $n>1$ or
    $-1\in N_{K/F}(K^\times)$.  Let $\delta\in K^\times$ be any
    exceptional element of $K/F$.

    Then the $\Fp[G]$-module $J$ decomposes as
    \begin{equation*}
        J=X\oplus Y, \quad Y = Y_n \oplus Y_{n-1} \oplus \dots
        \oplus Y_0,
    \end{equation*}
    where
    \begin{enumerate}
        \item $X$ is the cyclic $\Fp[G]$-module generated by
        $[\delta]$, with dimension $p^m+1$;
        \item $Y_i$ is a direct sum of cyclic $\Fp[G]$-modules of
        dimension $p^i$; and
        \item for $i \in \{0, \dots, n-1\}$
        \begin{equation*}
            [K_i^\times]=\left\{%
            \begin{array}{ll}
                X^{(\sigma-1)}\oplus Y^{H_i}, & m \le i; \\
                X^{(\sigma-1)(\sigma^{p^i}-1)^{p^{m-i}-1}}\oplus
                Y^{H_i}, & i<m. \\
            \end{array}%
            \right.
        \end{equation*}
    \end{enumerate}
    (Here $X^{(\sigma-1)}$ and $X^{(\sigma-1)
    (\sigma-1)^{p^{m-i}-1}}$ denote the images of $X$ under the
    action of $(\sigma-1)$ and
    $(\sigma-1)(\sigma^{p^i}-1)^{p^{m-i}-1}$, respectively.)
\end{theorem}

As before, let
\begin{equation*}
    e_i = \dim_{\Fp}
        \left(\left[N_{K_i/F}\left(K_i^\times\right)
    \right]/\left[N_{K_{i+1}/F}\left(K_{i+1}^\times
    \right)\right]\right), \quad 0\le i<n,
\end{equation*}
and let $e_n = \dim_{\Fp}[N_{K/F}(K^\times)]$.

\begin{corollary}\label{co:x}
    For each $m < i \le n$
    \begin{equation*}
        [N_{K_i/F}(K_i^\times)] = (Y_i+Y_{i+1}+\dots+Y_n)^G,
    \end{equation*}
    and, if $m\ge 0$, for each $0\le i\le m$
    \begin{equation*}
        [N_{K_i/F}(K_i^\times)] = (X+Y_i+Y_{i+1}+\dots+Y_n)^G.
    \end{equation*}

    For $i\neq m$
    \begin{equation*}
        \rank_{\Fp[G/H_i]} Y_i = e_i,
    \end{equation*}
    while if $m\ge 0$
    \begin{equation*}
        1 + \rank_{\Fp[G/H_m]} Y_{m} = e_{m}.
    \end{equation*}
\end{corollary}

Finally, we present several interesting conditions equivalent to
$m=i(K/F)$ being a particular element of the subset of field
indices $\mathcal{E} = \{-\infty,0,\dots,n-1\}$.  To express these
conditions, we define $-\infty \dotplus 1 = 0$ and, for $e \in
\mathcal{E}$ with $e \geq 0$, we define $e \dotplus 1 = e+1$.
We also set $N_{K_{n-1}/F}(K^\times_{-\infty})$ to be $\{1\}$.

\begin{theorem}\label{th:equivconds}
    Suppose that $\xi_p\in F$ and, if $p=2$, that either $n>1$ or
    $-1\in N_{K/F}(K^\times)$.  Then
    \begin{equation*}
    \begin{split}
        i(K/F) & = \min \left\{s \ \vert \ \xi_p \in
        N_{K/F}(K^\times) N_{K_{n-1}/F}(K_{s}^\times)\right\} \\
        &= \min\left\{ s \ \vert \ \xi_p \in N_{K/K_{s \dotplus
        1}} (K^\times) \right\} \\
        &= \min\left\{ s \ \vert \ \exists
        [\delta] \in J^{H_{s \dotplus 1}}, \ [N_{K/K_{s \dotplus 1}}
        \delta]_{K_{s \dotplus 1}} \neq [1]_{K_{s\dotplus 1}}
        \right\}.
    \end{split}
    \end{equation*}
\end{theorem}

From the equality
\begin{equation*}
    i(K/F) = \min \{ s \mid \xi_p \in N_{K/K_{s \dotplus
    1}}(K^\times) \}
\end{equation*}
above we obtain the following corollary.

\begin{corollary}\label{co:equivconds}
    Suppose that $\xi_p \in F$ and, if $p=2$, then either $n > 1$ or
    $-1\in N_{K/F}(K^\times)$. Then
    \begin{enumerate}
        \item for each $0 \le j \le n-1$
        \begin{equation*}
            i(K/K_j)=\left\{%
            \begin{array}{ll}
                i(K/F)-j, & 0\le j\le i(K/F) \\
                -\infty, & 0\le i(K/F)< j \\
                -\infty, & i(K/F)=-\infty \\
            \end{array}%
            \right.
        \end{equation*}
        \item for each $0<j\le n-1$, we have $i(K_j/F)=-\infty$.
    \end{enumerate}
\end{corollary}

One can also connect these equalities with the existence of
solutions of particular Galois embedding problems. This connection
will be pursued in a forthcoming paper. $X$ summands also lead
naturally to an investigation of the connections between $i(K/F)$
and the index of certain cyclotomic cyclic algebras, as well as the
behavior of $i(K/F)$ under adjoining roots of unity \cite{MSS1}. We
also show in \cite{MSS1} that all pairs $(i(K/F),n)$ with $n\in \N$
and $i(K/F)\in \{-\infty, 0, \dots, n-1\}$  are realizable for a
suitable field extension $K/F$ with $\xi_p\in F$. We plan a further
study of the behavior of $i(K/F)$ under base extension.

The proofs of Theorems~\ref{th:nox} and \ref{th:x} are inductive,
resting on the base case $n=1$ for Theorem~\ref{th:nox} and two base
cases $n=1$ and $p=2$, $n=2$ for Theorem~\ref{th:x}.  In these base
cases as well as the inductive proof, we employ lemmas which
establish the structure of the fixed submodule $J^G$ of $J$---in
particular, whether this fixed submodule is no more than the image
of the $p$th-power classes of the base field $F$---and specify which
of these elements are norms.

In fact, these lemmas reflect what has emerged, both in this work
as well as in the work on determining the entire Milnor $K$-theory
modulo $p$ as a Galois module (see \cite{BLMS} and \cite{LMS}), as
two essential foundational ingredients in the proof.  The first is
Hilbert's Theorem~90, which in our situation may be viewed as a
principle saying that we have enough norms. Indeed, Hilbert~90
tells us that the kernel of the norm map is as small as possible.
In order to use Hilbert~90 effectively, we need again and again
the technical refinements of this principle telling us that
certain elements in a group of $p$th-power classes are norms. In
this work these refinements, for example, begin with
Lemmas~\ref{le:fixed1}, \ref{le:fixed2}, and \ref{le:fixednew}
(identifying some fixed elements as norms), and are completed in
the full proofs of Theorems~\ref{th:nox} and \ref{th:x}.

The second essential ingredient is control of the image of
$p$th-power classes of the base field in the group of $p$th-power
classes of our field extension, which in this work is obtained
from Lemma~\ref{le:exact} (the Exact Sequence Lemma) and its
technical relative Lemma~\ref{le:fixed} (the Fixed Submodule
Lemma). In this paper, both of these principles are elementary,
but they are more sophisticated in the higher Milnor $K$-theory
case.  It is remarkable that one requires only repetitions of
these two principles in order to determine fully the Galois module
structure of the modules in question.  Drawing out the structure
from only these two first principles, however, does not come
without cost, and a number of technical observations turn out to
be necessary for us to fit the puzzle pieces together.

When $\xi_p\in F$, we need additional information to determine
when an element $[\gamma]\in J^{H_i}$ lies in $[K_i^\times]$ or is
instead an exceptional element.  We begin by standardizing choices
of the $a_i$ in the presentations of subfields $K_{i+1} =
K_i(\sqrt{a_i})$ in section~\ref{se:aipres}.  Then, in
section~\ref{se:submods}, we collect several results used in
identifying elements of $[K_i^\times]$.  These are the
Submodule-Subfield Lemma~\eqref{le:sub} for free components, the
Norm Lemma~\eqref{le:norm} for comparisons among norms from $K$ to
various $K_i$ (in order to determine when an exceptional element
for $K/F$ is an exceptional element for $K/K_i$), and the Proper
Subfield Lemma~\eqref{le:prop} for elements that generate
sufficiently small cyclic submodules.

In section~\ref{se:reptheory}, we present lemmas which we use to
manipulate $\Fp[G]$-representa\-tions formally: the Inclusion
Lemma~\eqref{le:incl}, the Exclusion Lem\-ma~\eqref{le:excl}, and
the Free Complement Lemma~\eqref{le:free}.

We begin the proof by proving the base cases for an induction in
section~\ref{se:base}. Our inductive strategy is first to show that
$J$ contains a sufficiently large direct sum of $\Fp[G]$-submodules
of $p$th-power dimensions. We do so in section~\ref{se:nox} in
Proposition~\ref{pr:y}, the result of which is already enough to
prove Theorem~\ref{th:nox}. When $\xi_p\in F$ and, if $p=2$, $n>1$,
we also need to establish the dimension of the $X$ component and
connect notions of exceptional elements for subextensions $K/K_i$.
We do so in section~\ref{se:exc}. In section~\ref{se:x}, we prove an
analogue of Proposition~\ref{pr:y} which establishes
Theorem~\ref{th:x} without the independence of $X$ and $Y$, and then
we prove Theorem~\ref{th:x} fully.  Finally, in
section~\ref{se:equivconds}, we prove Theorem~\ref{th:equivconds}.

For the reader's convenience, we have made our paper
self-contained; in particular, it is independent from \cite{MS}.

\tableofcontents

\section{Notation and Lemmas}

\subsection{$\Fp[G]$-modules}\label{se:reptheory}

Let $G$ be a cyclic group of order $p^n$ with generator $\sigma$.
For an $\Fp[G]$-module $U$, let $U^G$ denote the submodule of $U$
fixed by $G$, and for an arbitrary element $u\in U$, let $l(u)$
denote the dimension of the $\Fp[G]$-submodule of $U$ generated by
$u$. Denote by $N$ the operator $(\sigma-1)^{p^n-1}$ acting on $U$.
For an $\Fp[G]$-module $V$ and an element $\gamma \in V$, let
$\langle \gamma \rangle$ denote the $\Fp$-subspace of $V$ spanned by
$\gamma$, and let $M_\gamma$ denote the cyclic $\Fp[G]$-module
generated by $\gamma$. If $[\gamma]$ is an element of
$K^\times/K^{\times p}$ represented by $\gamma\in K^\times$, we
write $M_\gamma$ instead of $M_{[\gamma]}$.  Similarly we also write
$l(\gamma)$ instead of $l([\gamma])$.

We will usually use additive notation for general
$\Fp[G]$-modules, switching to multiplicative notation when
considering the specific module $J=K^\times/K^{\times p}$.
However, occasionally even in this case we employ additive
notation, in particular writing $\{0\}$ to denote $\{[1]\}$.

\begin{lemma}[Inclusion Lemma]\label{le:incl}
    Let $U$ and $V$ be $\Fp[G]$-modules contained in an
    $\Fp[G]$-module $W$.  Suppose that $(U+V)^G\subset U$ and for
    each $w \in (U+V)\setminus (U+V)^G$ there exists $u\in U$ such
    that
    \begin{equation*}
        (\sigma-1)^{l(w)-1}(w)=N(u).
    \end{equation*}
    Then $V\subset U$.
\end{lemma}

\begin{proof}
    Let $\{T_i\}_{i=1}^s$ be the socle series of $U+V$:
    $T_1=(U+V)^G$ and $T_{i+1}/T_i=((U+V)/T_i)^G$, and let $s$ be
    the least natural number such that $T_s=U+V$. Observe that
    since $(\sigma-1)^{p^n}=0$, we have $s\le p^n$.  We prove the
    lemma by induction on the socle series.

    By hypothesis, $T_1\subset U$.  Assume now that $T_i\subset U$
    for some $i<s$. Then for each $w\in T_{i+1}\setminus T_i$ we
    have $l(w)= i+1$ and $(\sigma-1)^{l(w)-1}(w) = N(u) =
    (\sigma-1)^{p^n-1}(u)$ for some $u\in U$.  Therefore
    \begin{equation*}
        (\sigma-1)^{l(w)-1} \big(w-(\sigma-1)^{p^n-l(w)}(u)\big) =
        0.
    \end{equation*}
    Therefore $w-(\sigma-1)^{p^n-l(w)}(u)\in T_i\subset U$.  Hence
    $w\in U$ and $T_{i+1}\subset U$. Therefore $U+V=U$ and $V
    \subset U$ as required.
\end{proof}

\begin{lemma}[Exclusion Lemma]\label{le:excl}
    Let $U$ and $V$ be $\Fp[G]$-modules contained in an
    $\Fp[G]$-module $W$.  Suppose that $U^G\cap V^G=\{0\}$.  Then
    $U+V=U\oplus V$.
\end{lemma}

\begin{proof}
    Let $Z=U\cap V$ and suppose that $y\in Z\setminus \{0\}$.  Let
    \begin{equation*}
        z=(\sigma-1)^{l(y)-1}(y)\neq 0.
    \end{equation*}
    Then $z\in U^G\cap V^G$, a contradiction. Hence $U\cap
    V=\{0\}$ and $U+V=U\oplus V$.
\end{proof}

We will use the Exclusion Lemma~\eqref{le:excl} in the case of
several or even infinitely many modules, as well, since by a simple
argument one can always reduce such a case to the case of two
modules. We indicate this argument in the proof of the following
lemma and omit it later on.

The following lemma follows from the fact that each free
$\Fp[G]$-module is injective.  (See \cite[Theorem~11.2]{Ca}.) We
shall, however, provide a direct proof.

\begin{lemma}[Free Complement Lemma]\label{le:free}
    Let $V\subset U$ be free $\Fp[G]$-modules.  Then there exists
    a free $\Fp[G]$-submodule $\tilde V$ of $U$ such that $V\oplus
    \tilde V=U$.
\end{lemma}

\begin{proof}
    Let $Z$ be a complement of $V^G$ in $U^G$ as $\Fp$-vector
    spaces, and let $\Zc$ be an $\Fp$-base of $Z$.  For each $z\in
    \Zc$, there exists $u(z)$ such that $z=N(u(z))$.  Let $M(z)$
    be the $\Fp[G]$-submodule of $U$ generated by $u(z)$. Then
    $M(z)$ is a free $\Fp[G]$-submodule.  Moreover, its fixed
    submodule $M(z)^G$ is the $\Fp$-vector subspace generated by
    $z$.

    We claim that the $M(z)$, $z\in \Zc$, are independent.
    First we show by induction on the number of modules that a
    finite set of modules $M(z)$ is independent.  The base case is
    trivial.  Now let $W=M(z)\cap \sum_{z'\neq z} M(z')$.  Now
    by the inductive assumption on independence, $(\sum_{z'\neq
    z} M(z'))^G = \sum_{z'\neq z} M(z')^G$, and for each
    $z$, $M(z)^G = \langle z\rangle$. Since the $z$ form an
    $\Fp$-base for $Z$, we obtain $W^G=\{0\}$.  The
    Exclusion Lemma~\eqref{le:excl} then gives
    $M(z)+\sum_{z'\neq z} M(z') = M(z) \oplus \sum_{z'\neq
    z} M(z')$.

    The case of an infinite sum follows from the same argument,
    since the fact that $m\in M(z)^G \cap \sum_{z'\neq
    z} M(z')^G$ forces $m$ to be a finite sum of elements
    $m(z')$.  Hence the $M(z)$, $z\in \Zc$, are independent.

    Set $\tilde V := \oplus_{z\in \Zc} M(z)$.  Then $\tilde V$ is
    a free $\Fp[G]$-submodule of $U$ and $\tilde V^G=Z$.  By the
    Exclusion Lemma~\eqref{le:excl}, we have $V+\tilde V =
    V\oplus \tilde V$ and $(V\oplus \tilde V)^G=V^G\oplus \tilde
    V^G = U^G$.

    Now let $u\in U$ be arbitrary and let $M$ be the cyclic
    $\Fp[G]$-submodule of $U$ generated by $u$.  Then $(M+V+\tilde
    V)^G \subset U^G \subset V+\tilde V$.  Moreover, for any $m\in
    (M+V+\tilde V) \setminus (M+V+\tilde V)^G$
    \begin{equation*}
        (\sigma-1)^{l(m)-1}(m)\in (M+V+\tilde V)^G \subset U^G =
        (V+\tilde V)^G = N(V + \tilde V)
    \end{equation*}
    by the freeness of $V$ and $\tilde V$. By the Inclusion
    Lemma~\eqref{le:incl}, then, $M\subset V + \tilde V$. Hence
    $U=V\oplus \tilde V$.
\end{proof}

\subsection{Kummer Subfields of $K/F$ and Exceptional
Elements}\label{se:aipres}

Suppose that $\xi_p\in F$.  In this case we have Kummer theory and
may organize presentations of the extensions $K_{i+1}/K_i$ as
follows.

\begin{proposition}[Subfield Generators]\label{pr:subgen}
    We may choose $a_i\in K_i^\times$, $0\le i < n$ such that
    \begin{itemize}
        \item $K_{i+1} = K_i(\root{p}\of{a_i})$ and
        \item $N_{K_i/K_j} a_i=a_j$ for all $0\le j < i < n$.
    \end{itemize}
\end{proposition}

In what follows we will assume that the choices of $a_i$ have been
made according to Proposition~\ref{pr:subgen}, and we set $a=a_0$.

We prove this result by means of the following

\begin{lemma}\label{le:p2norms}
    Suppose that $\xi_p\in K$ and let $L'/K$ be a cyclic extension
    of degree $p^2$ with $L/K$ the intermediate extension of
    degree $p$. Then, for every $b\in L$ with
    $L'=L(\root{p}\of{b})$, we have $L=K(\root{p}\of{
    N_{L/K}(b)})$.
\end{lemma}

\begin{proof}
    Let $\sigma$ be a generator of $\Gal(L'/K)$. For each $i \in\{1,
    2, \dots, p-1\}$, we have
    \begin{equation*}
        \left(\root{p}\of{b}\right)^{\sigma^i} =
        \root{p}\of{b^{\sigma^i}}
    \end{equation*}
    for a suitable choice of a $p$th root of $b^{\sigma^i}$. Hence
    \begin{equation*}
        \left(\root{p}\of{b}\right)^{1+\sigma+ \dots+\sigma^{p-1}}
        = \root{p}\of{b^{1+\sigma+\dots+\sigma^{p-1}}} =
        \root{p}\of{N_{L/K}(b)}\in L'
    \end{equation*}
    for a suitable choice of a $p$th root of $N_{L/K}(b)$.

    Observe that since $\xi_p\in K$ the equality
    \begin{equation*}
        \left(\root{p}\of{b}\right)^{(1+\sigma+\dots+\sigma^{p-1})
        (\sigma-1)} = \root{p}\of{b}^{\sigma^p-1} =
        \root{p}\of{N_{L/K}(b)}^{\sigma-1}
    \end{equation*}
    is independent of the choice of $p$th roots. Moreover, since
    $L'=L(\root{p}\of{b})$ and $\sigma^p$ generates $\Gal(L'/L)$,
    we see that $\root{p}\of{b}^{\sigma^p-1} \neq 1$. Hence we
    conclude that $L=K(\root{p}\of{N_ {L/K}(b)})$.
\end{proof}

\begin{proof}[of Proposition~{\rm\ref{pr:subgen}}]
    By Kummer theory, there exists $a_{n-1}\in K_{n-1}^\times$
    such that $K_n = K_{n-1} (\root{p}\of{a_{n-1}})$. Then
    inductively define
    \begin{equation*}
        a_{n-i}=N_{K_{n-i+1}/K_{n-i}}(a_{n-i+1})
    \end{equation*}
    for $i\in\{2, \dots, n\}$. Applying the lemma to extensions
    $K_{n-i+2}/K_{n-i}$, we have the results.
\end{proof}

Our definition of exceptional elements makes use of a subset of
the set  $\{\delta\in K^\times \ \vert\ [N_{K/F}(\delta)]_F \neq
[1]_F\}$. In general, however, this latter set may be empty.
Consider, for example, the extension $\C/\R$, for which
$N_{\C/\R}(\C^\times)\subset \R^{\times 2}$.  The next
proposition shows that under the conditions we require in the
definition of exceptional elements, this set is never empty and
therefore exceptional elements exist.

\begin{proposition}\label{pr:excexist}
    Let $\xi_p\in F$ and, if $p=2$, that $n>1$ or $-1\in
    N_{K/F}(K^\times)$. Then an exceptional element $\delta$ exists.
\end{proposition}

\begin{proof}
    Consider $\delta = \root{p}\of{a_{n-1}}$. If $p>2$ then
    $N_{K/K_{n-1}}(\delta) = a_{n-1}$ and hence $N_{K/F}(\delta) =
    a_0 = a$. Now if $p=2$ then $N_{K/K_{n-1}}(\delta) = -a_{n-1}$
    and for $n>1$ we similarly have $N_{K/F}(\delta) = a_0 = a$. If
    $p=2$ and $n=1$, then $-a = N_{K/F}(\sqrt{a})$ and hence $-1\in
    N_{K/F}(K^\times)$ if and only if $a\in N_{K/F}(K^\times)$.
    Consequently, under our hypothesis, exceptional elements always
    exist.
\end{proof}

\subsection{The Fixed Submodule $J^G$ of $J$}\label{se:fixed}

Recall that we write $[F^\times]$ for $F^\times K^{\times
p}/K^{\times p}\subset J$.

The following lemmas generalize \cite[Lemma 2 and Remark 2]{MS}:

\begin{lemma}[Fixed Submodule Lemma]\label{le:fixed}\
\begin{enumerate}
    \item If $\xi_p\not\in N_{K/F}(K^\times)$
    \begin{equation*}
        J^G = [F^\times].
    \end{equation*}
    \item If $\xi_p\in N_{K/F}(K^\times)$
    \begin{equation*}
        J^G = \langle [\delta]\rangle \oplus [F^\times],
    \end{equation*}
    where $\delta\in K^\times$ with $\delta^{\sigma-1}=\lambda^p$,
    $N_{K/F}(\lambda)$ is a primitive $p$th root of unity, and
    $[N_{K/F}(\delta)]_F=[a]_F$.  In particular, $\delta$ is
    an exceptional element of $K/F$ and $i(K/F)=-\infty$.
\end{enumerate}
\end{lemma}

\begin{proof}
    Suppose that $\theta\in K^\times$ such that $[\theta]\in
    J^G$. Then $\theta^{\sigma-1}=\lambda^p$ for some
    $\lambda\in K^\times$, and hence $N_{K/F}(\lambda)^p=1$.
    Therefore $N_{K/F}(\lambda)$ is a $p$th root of unity.

    Now consider the first case, $\xi_p\not\in N_{K/F}(K^\times)$.
    Then $N_{K/F}(\lambda)=1$, because otherwise $\xi_p$ would be
    the norm of a suitable power of $\lambda$.  From Hilbert~90 we
    see that $\theta^{\sigma-1}=(k^p)^{\sigma-1}$ for some $k\in
    K^\times$.  We conclude that $\theta/k^p\in F^\times$ and hence
    $[\theta]=[f]$ for some $f\in F^\times$.  Therefore if
    $\xi_p\not \in N_{K/F}(K^\times)$ then $J^G=[F^\times]$ as
    required.

    Now assume that $\xi_p\in N_{K/F}(K^\times)$.  Then $\xi_p =
    N_{K/F}(\lambda)$ for some $\lambda\in K^\times$ and by
    Hilbert~90 there exists an element $\delta\in K^\times$ such
    that $\delta^{\sigma-1}=\lambda^p$.  Then the $\Fp[G]$-submodule
    of $J$ generated by $[\delta]$ and $\epsilon(F^\times)$ is
    isomorphic to $[F^\times]\oplus \langle [\delta]\rangle$.

    By \cite[Theorem 3]{A}, $K(\root{p}\of{\delta})$ is a cyclic
    extension of $F$ of degree $p^{n+1}$.  Then repeated
    application of Lemma~\ref{le:p2norms} gives
    \begin{equation*}
        K_{n-i}=K_{n-i-1}\left(\root{p}\of{ N_{K_{n}/K_{n-i-1}}
        (\delta)}\right)
    \end{equation*}
    for $i\in\{0, 1, \dots, n-1\}$.  Hence we obtain $K_1=F(\root{p}
    \of {N_{K/F} (\delta)})$.  By Kummer theory, $\langle
    [N_{K/F}(\delta)]_F \rangle = \langle [a]_F\rangle$ as subgroups
    of $F^\times/ {F^{\times p}}$.  By replacing $\delta$ with
    another power if necessary, then, $[N_{K/F}(\delta)]_F=[a]_F$
    and $\delta^{\sigma-1}=\lambda^p$, where $N_{K/F}(\lambda)$ is a
    primitive $p$th root of unity.  We have $[\delta]^{(\sigma-1)} =
    [1]$ and so by definition $\delta$ is exceptional for $K/F$.
    Moreover we see from the definition of $i(K/F)$ that
    $i(K/F)=-\infty$ in this case.

    Now for each $[\theta]\in J^G$, $\theta^{\sigma-1}=\nu^p$
    with $N_{K/F}(\nu)=N_{K/F}(\lambda)^c$ for some $c\in \Z$.
    Then we have $(\theta\delta^{-c})^{\sigma-1}= \nu^p
    \lambda^{-pc}$. Because $N(\nu\lambda^{-c}) =1$, from
    Hilbert~90 we see that there exists $\omega\in K^\times$ such
    that $\omega^{\sigma-1} = \nu \lambda^{-c}$. Hence
    $(\theta\delta^{-c})^{\sigma-1} = (\omega^p)^{\sigma-1}$ and
    we see that $[\theta] \in [F^\times] + [\delta]^c$.  Therefore
    $J^G\cong [F^\times] \oplus \langle [\delta]\rangle$, as
    required.
\end{proof}

\begin{lemma}[Exact Sequence Lemma]\label{le:exact}
    There is an exact sequence
    \begin{equation*}
        1\to A\to F^\times/F^{\times p} \xrightarrow{\epsilon}
        J^G \xrightarrow{N} A
    \end{equation*}
    where $A=(F^\times \cap K^{\times p})/F^{\times p}$, $\epsilon$
    is the natural homomorphism induced by the inclusion $F^\times
    \to K^\times$, and $N$ is the homomorphism induced by the norm
    map $N_{K/F}\colon K^\times\to F^\times$.

    \begin{itemize}
        \item If $\xi_p\not\in F$, then $A=1$.
        \item If $\xi_p\in F$, $A=\langle [a]\rangle$, in which
        case the map $N$ is surjective if and only if $\xi_p\in
        N_{K/F}(K^\times)$.
    \end{itemize}
\end{lemma}

\begin{remark*}
    For future reference it is convenient to translate exactness
    at $J^G$ in the sequence above, as follows:
    \begin{enumerate}
        \item[(i)] If $[\beta]\in [F^\times]$ then
        $[N_{K/F}(\beta)]_F =[1]_F$.
        \item[(ii)] If $[\beta]\in J^G$ and $[N_{K/F}(\beta)]_F =
        [1]_F$ then $[\beta]\in [F^\times]$.
    \end{enumerate}
\end{remark*}

\begin{proof}
    If $\xi_p\in F$, then Kummer theory implies that the first
    occurrence of $A$ in the exact sequence above is equal to
    $A=\langle[a]_F\rangle$.  Otherwise, suppose that $\xi_p
    \notin F$. If $\chr(F)=p$ then no primitive $p$th root of
    unity lies in the algebraic closure of $F$, whence
    $\xi_p\notin K$.  If $\chr(F)\neq p$, then since
    $2\le[F(\xi_p):F]\leq p-1$ and $[K:F]=p^n$, we similarly obtain
    $\xi_{p}\notin K$. In any case, then, $\xi_p\notin K$. Assume
    that $k^p = f \in F^\times$.  Then $(k^p)^{\sigma-1} =
    (k^{\sigma-1})^p = 1$, whence $k^{\sigma-1}$ is a $p$th root
    of unity, which must be $1$. Hence $k^{\sigma-1}=1$, and we
    deduce $k\in F$ and $f\in F^{\times p}$.  Therefore $A=1$.

    The Fixed Submodule Lemma~\eqref{le:fixed} then gives exactness
    at $J^G$ and that $N$ is surjective if and only if either
    $\xi_p\not\in F$ or $\xi_p\in N_{K/F}(K^\times)$. Exactness at
    $F^\times/F^{\times p}$ follows from Kummer theory.
\end{proof}

\subsection{$\Fp[G]$-Submodules of $J$}\label{se:submods}

\begin{lemma}[Submodule-Subfield Lemma]\label{le:sub}
    Let $U$ be a free $\Fp[G]$-sub\-module of $J$ and
    $i\in\{0,1,\dots,n\}$. Then
    \begin{equation*}
        U^{H_i} = U^{(\sigma-1)^{p^n-p^i}} =
        U\cap [N_{K_n/K_i}K_n^\times] =
            U\cap [K_i^\times].
    \end{equation*}
\end{lemma}

\begin{proof}
    Suppose $[u]\in U^{H_i}$.  Then
    $[u]^{(\sigma^{p^i}-1)}=[u]^{(\sigma-1)^{p^i}}=[1]$, so
    $l(u)\le p^i$. Since $U$ is free, $[u]=[\tilde
    u]^{(\sigma-1)^{p^n-l(u)}}$ for some $[\tilde u]\in U$.  In
    particular
    \begin{equation*}
        [u]=([\tilde u]^{(\sigma-1)^{p^i-l(u)}})^
        {(\sigma-1)^{p^n-p^i}}.
    \end{equation*}
    Hence $U^{H_i}\subset U^{(\sigma-1)^{p^n-p^i}}$.  Now suppose
    that $[u]=[\tilde u]^{(\sigma-1)^{p^n-p^i}}$.  Then since
    \begin{equation*}
        [\tilde u]^{(\sigma-1)^{p^n-p^i}} = [N_{K_n/K_i}({\tilde
        u})],
    \end{equation*}
    $U^{(\sigma-1)^{p^n-p^i}}\subset U\cap [N_{K_n/K_i}
    K_n^\times]\subset U\cap [K_i^\times]$.

    Finally suppose that $[u]\in U\cap [K_i^\times]$.  Then
    $[u]\in U^{H_i}$ and we see that all of our inclusions
    above are actually equalities.
\end{proof}

\begin{remark*}
    If $U$ is a free $\Fp[G]$-module, then $U$ is also a free
    $\Fp[H_i]$-module.  But then $H^2(H_i,U)=\{0\}$.  Hence
    $U^{H_i}=N_i(U):=$ the image of the norm operator $N_i$.  Thus
    $U^{H_i}=U^{(\sigma-1)^{p^n-p^i}}$ as required.
\end{remark*}

Just as with $F=K_0$, denote elements of the $\Fp[G/H_i]$-module
$J(K_i) = K_i^\times/K_i^{\times p}$ by $[\gamma]_{K_i}$, $\gamma
\in K_i^\times$.

\begin{lemma}[Norm Lemma]\label{le:norm}
    For all elements $[\gamma]\in J$ with $l(\gamma) < p^n$,
    we have $[N_{K/F}(\gamma)]_F \in \langle [a]_F \rangle$.

    Now suppose additionally that $l(\gamma)\le p^n-p^i$ for some
    $0\le i< n$. Then $[N_{K/K_i}(\gamma)]_{K_i} \in \langle
    [a_i]_{K_i}\rangle$, and $[N_{K/F}(\gamma)]_F=[a]_F^s$ if and
    only if $[N_{K/K_i}(\gamma)]_{K_i}=[a_i]_{K_i}^s$.
\end{lemma}

\begin{proof}
    For the first statement, observe that $(1+\sigma+\dots +
    \sigma^{p^n-1}) \equiv (\sigma-1)^{p^n-1}$ on $J$, and hence
    $[N_{K/F}(\gamma)] = [\gamma]^{(\sigma-1)^{p^n-1}}$.  Since
    $l(\gamma) < p^n$, $[N_{K/F}(\gamma)] = [1]$.  Therefore
    $N_{K/F}(\gamma) \in F^\times \cap K^{\times p}$, which by
    Kummer theory is the union $\cup_{j=0}^{p-1} a^jF^{\times p}$.
    We obtain $[N_{K/F}(\gamma)]_F \in \langle [a]_F \rangle$.

    For the second statement, observe first that if $[\gamma]=[1]$
    then the lemma is trivial. Otherwise, consider $J$ as an
    $\Fp[H_i]$-module and let $\tau = \sigma^{p^i}$.  The
    $\Fp[H_i]$-module generated by $[\gamma]$ has dimension equal
    to $t$, where $[\gamma]^{(\tau-1)^t}=[1]$ and
    $[\gamma]^{(\tau-1)^{t-1}}\neq [1]$. Since $\tau-1\equiv
    (\sigma-1)^{p^i}$ on $J$, this condition is equivalent to
    $(t-1)p^i < l(\gamma) \le tp^i$.  Since $l(\gamma) \le
    (p^{n-i}-1)p^i$, the dimension $t$ is strictly less than
    $p^{n-i}$. Hence $(\tau-1)^{p^{n-i}-1}$ annihilates the cyclic
    $\Fp[H_i]$-module generated by $\gamma$, and so its length, as
    an $\Fp[H_i]$-module, is less than $p^{n-i}$.

    Applying the first statement in the case of the cyclic
    extension $K/K_i$, we have $[N_{K/K_i}(\gamma)]_{K_i}\in
    \langle [a_i]_{K_i} \rangle$. Now because $N_{K_i/F}(a_i) = a$
    and
    \begin{equation*}
        [N_{K/F}(\gamma)]_F = N_{K_i/F}
        ([N_{K/K_i}(\gamma)]_{K_i}),
    \end{equation*}
    we have $[N_{K/K_i}(\gamma)]_{K_i} = [a_i]_{K_i}^s$ if and
    only if $[N_{K/F}(\gamma)]_F = [a]_F^s$.
\end{proof}

\begin{remark*}
    Occasionally, we will cite the Norm Lemma~\eqref{le:norm} as
    an abbreviation of the simple argument, at the end of the
    lemma's proof, which shows that
    \begin{equation*}
        [N_{K/F}(\gamma)]_F = [a]_F^s \text{ if and only if }
        [N_{K/K_i}(\gamma)]_{K_i} = [a_i]_{K_i}^s.
    \end{equation*}
\end{remark*}

\begin{lemma}[Proper Subfield Lemma]\label{le:prop}
    Let $[z]\in J^{H_i}$, $i<n$.  Then $[z]\in [K_i^\times]$ if
    and only if $[N_{K/F}(z)]_F=[1]_F$.
\end{lemma}

\begin{proof}
    If $[z]\in J^{H_i}$, then $[z]^{(\sigma^{p^i}-1)}=[1]$.  Since
    $(\sigma^{p^i}-1)\equiv (\sigma-1)^{p^i}$ on $J$, $l(z)\le
    p^i$.  Further assume that $[N_{K/F}(z)]_F = [1]_F$.

    Consider $J$ as an $\Fp[H_i]$-module. Then from the Exact Sequence
    Lemma~\eqref{le:exact} applied to the field
    extension $K/K_i$, we see that
    \begin{equation*}
        [z]\in [K_i^\times] \text{\ \ \ or\ \ \ } [z]\in \left(\langle
    [\delta]\rangle\oplus [K_i^\times]\right) \setminus [K_i^\times]
    \end{equation*}
    according to whether
    \begin{equation*}
        [N_{K/K_i}(z)]_{K_i} = [1]_{K_i} \text{\ \ \ or\ \ \ }
        [N_{K/K_i}(z)]_{K_i} \neq [1]_{K_i}.
    \end{equation*}
    (Here $\delta\in K^\times$ with $\delta^{\sigma-1}
    =\lambda^p$, $N_{K/K_i}(\lambda)$ is a primitive $p$th root of
    unity, and $[N_{K/K_i}(\delta)]_{K_i} = [a_i]_{K_i}$.)

    Therefore if $[z]\notin [K_i^\times]$ then $[N_{K/K_i}(z)]_{K_i}
    = [a_i]_{K_i}^c$ for $c\not\equiv 0 \bmod p$, and by the Norm
    Lemma~\eqref{le:norm}, $[N_{K/F}(z)]_F=[a]_F^c \neq [1]_F$ --- a
    contradiction.

    Conversely, if $[z]\in [K_i^\times]$ then $[z]\in J^{H_i}$ and
    \begin{equation*}
        [N_{K/F}(z)]_F = [N_{K_i/F}(z)]_F^{p^{n-i}} = [1]_F,
    \end{equation*}
    since $n>i$.
\end{proof}

\subsection{Fixed Submodules of Cyclic Submodules of $J$}
\label{se:norms}

\begin{lemma}[First Fixed Elements are Norms Lemma]\label{le:fixed1}
    Suppose that $p>2$, $n=1$, $[\gamma]\in J$, $2\le l(\gamma) <
    p$, and that one of the following holds:
    \begin{itemize}
        \item $\xi_p\not\in F$
        \item $\xi_p\in F$ and $l(\gamma) \ge 3$
        \item $\xi_p\in F$, $l(\gamma)=2$, and $\gamma$ is
        unexceptional.
    \end{itemize}
    Then there exists $[\alpha] \in J$ such that $M_\gamma^G =
    \langle N[\alpha]\rangle$.
\end{lemma}

\begin{proof}
    First suppose that $\xi_p\not\in F$.  We show by induction on $i
    \in \{l(\gamma), \dots, p\}$ that there exists an element
    $\alpha_i\in K^\times$ with $\langle
    [\alpha_i]^{(\sigma-1)^{i-1}}\rangle = M_\gamma^G$. Then since
    $(\sigma-1)^{p-1} \equiv 1+\sigma+\dots+\sigma^{p-1}$ we may set
    $\alpha := \alpha_p$ and the proof of the first item will be
    complete.  If $i=l(\gamma)$ we set $a_i = \gamma$. Assume now
    that $l(\gamma)\le i< p$ and that our statement is true for $i$.

    Set $c=N_{K/F}(\alpha_i)$.  Since $[\alpha_i]^{
    (\sigma-1)^{p-1}} = [c]$ and $i<p$, we see that $[c]=[1]$.
    Then $c\in F^\times\cap K^{\times p}$, which by the Exact
    Sequence Lemma~\eqref{le:exact} is equal to $F^{\times p}$.
    Hence $c=f^p$ for some $f\in F^\times$. Then
    $N_{K/F}(\alpha_i/f)=1$. By Hilbert~90 there exists an element
    $\omega\in K^\times$ such that $\omega^{\sigma-1}=\alpha_i/f$.
    Then $\omega^{(\sigma-1)^2} = \alpha_i^{(\sigma-1)}$.  Since
    $l(\alpha_i)\ge 2$ and $\langle [\alpha_i]^{(\sigma-1)^{i-1}}
    \rangle = M_\gamma^G$, $\langle [\omega]^{(\sigma-1)^{i}}
    \rangle = M_\gamma^G$ and we may set $\alpha_{i+1}=\omega$.
    Our induction is complete.

    Now suppose that $\xi_p\in F$, $K=F(\root{p}\of{a})$, and
    $\root{p}\of{a}^{\sigma}=\xi_p\root{p}\of{a}$.  First assume
    that $l(\gamma)\ge 3$. As before, we show by induction on $i$
    that there exists an element $\alpha_{i}\in K^\times$ such that
    $\langle[\alpha_{i}]^{(\sigma-1)^{i-1}}\rangle =
    M_{\gamma}^{G}$. If $i=l(\gamma)$ we set $\alpha_{i}=\gamma$.
    Assume now that $l(\gamma)\leq i< p$ and that our statement is
    true for $i$.

    By the Norm Lemma~\eqref{le:norm} we have $[N_{K/F}(\alpha_i)]_F
    \in \langle [a]_F \rangle$. Hence $c := N_{K/F}(\alpha_i)
    =a^sf^p$ for some $f\in F^\times$ and $s\in \Z$. Then
    $N_{K/F}(\alpha_{i}/f \delta^s)=1$, where $\delta =
    \root{p}\of{a}$. By Hilbert~90 there exists an element
    $\omega\in K^\times$ such that $\omega^{\sigma-1}=\alpha_{i}/f
    \delta^s$. Then
    $\omega^{(\sigma-1)^{2}}=\alpha_{i}^{(\sigma-1)}/\xi_{p}^{s}$.
    Since $i \geq 3$, $\langle[\omega]^{(\sigma-1)^{i}}\rangle=
    \langle [\alpha_{i}]^{(\sigma-1)^{i-1}} \rangle =
    M_{\gamma}^{G}$ and we can set $\alpha_{i+1}:=\omega$.

    Assume then that $l(\gamma)=2$ and $\gamma$ is an
    unexceptional element of $K/F$. By the Norm
    Lemma~\eqref{le:norm}, $[N_{K/F}(\gamma)]_F \in \langle [a]_F
    \rangle$, and as before $c := N_{K/F}(\gamma) = a^sf^p$ for
    some $f\in F^\times$ and $s\in \Z$.

    Since $\gamma$ is unexceptional, either $s \equiv 0 \bmod p$,
    in which case $c = f^p$ for some $f \in F^{\times}$, or
    $[\gamma]^{\sigma-1} \notin [K_m^\times]$.  In the former
    case, $N_{K/F}(\gamma/f)=1$. By Hilbert~90 there exists an
    element $\omega\in K^\times$ such that $\omega^{\sigma-1} =
    \gamma/f$ and $\omega^{(\sigma-1)^2} = \gamma^{\sigma-1}$.
    Hence $\langle [\omega]^{(\sigma-1)^2}\rangle = M_\gamma^G$
    and we may invoke the statement for $\omega$ since
    $l(\omega)=3$.

    In the latter case, since $N_{K/F}(\gamma^{\sigma-1})=1$ and
    $[\gamma]^{\sigma-1}\in J^G$, from the Exact Sequence
    Lemma~\eqref{le:exact}(ii) we see that $[\gamma]^{\sigma-1}\in
    [F^\times]=[K_0^\times]$.  Hence $m<0$ so that $m=-\infty$.
    Thus there exists an element $\delta\in K^\times$ such that
    $[N_{K/F}(\delta)]_F \neq [1]_F$ and $[\delta]^{\sigma-1} =
    [1]$. Again using the Exact Sequence Lemma~\eqref{le:exact} we
    see that we may assume that $[N_{K/F}(\delta)]_F=[a]_F$ and
    $[\delta]^{\sigma-1}=[1]$.

    Now let $N_{K/F}(\delta) = ag^p$ for some $g\in F^\times$ and
    note $N_{K/F}(\gamma g^s/f\delta^s) = 1$. Then as before we
    have $\omega^{\sigma-1}=\gamma g^s/f \delta^s$ and $[\omega]^{
    (\sigma-1)^2} = [\gamma]^{(\sigma-1)} \neq [1]$. Hence
    $\langle [\omega]^{(\sigma-1)^2}\rangle = M_\gamma^G$ and we
    may invoke the statement for $\omega$ since $l(\omega)=3$.
\end{proof}

\begin{lemma}[Fixed Elements of Length 3 Submodules are
Norms Lem\-ma]\label{le:fixed2}
    Suppose that $p=2$, $n=2$, $[\gamma]\in J$, $l(\gamma) = 3$,
    and $[N_{K/F}(\gamma)]_F = [1]_F$. Then there exists $[\alpha]
    \in J$ such that $M_\gamma^G = \langle N[\alpha]\rangle$.
\end{lemma}

\begin{proof}
    Let $\beta=\gamma^{\sigma-1}$.  Then $l(\beta)=2$ and, since
    $\beta$ is in the image of $\sigma-1$, we have
    $[N_{K/F}(\beta)]_F = [1]_F$. Because $l(\beta)=2$ and
    $N_{K/K_1}$ is equivalent to $1+\sigma^2 \equiv (\sigma-1)^2$ on
    $J$, we see that $[N_{K/K_1}(\beta)]=[1]$ in $J$.  From the Norm
    Lemma~\eqref{le:norm} we conclude that $[N_{K/K_1}(\beta)]_{K_1}
    = [1]_{K_1}$, and by the Exact Sequence
    Lemma~\eqref{le:exact}(ii) applied to the $\Ft[H_{1}]$-module
    $J$, we see that $[\beta]\in [K_1^\times]$. Let $b\in
    K_1^\times$ such that $[b]=[\gamma]^{\sigma-1}$.

    Now set $c:=N_{K_1/F}(b)$. Observe that $[c] = [b]^{1+\sigma} =
    [\gamma]^{\sigma^2-1} = [N_{K/K_1}(\gamma)]$ and $\langle
    [c]\rangle \subset M_\gamma^G$. Hence $N_{K/K_1}(\gamma) = ck^2$
    for some $k\in K^\times$, and $k^2\in K_1^\times \cap K^{\times
    2}$.  By Kummer theory $k^2=a_1^sg^2$ for some $s\in \Z$ and
    $g\in K_1^\times$, whence $N_{K/K_1}(\gamma)= ca_1^sg^2$ and
    $[N_{K/F}(\gamma)]_F = [a]^s_F$. By hypothesis $s\equiv 0\bmod
    2$. Therefore $N_{K/K_1}(\gamma)=ch^2$ for some $h\in
    K_1^\times$.

    Now $N_{K/F}(\gamma)=N_{K_1/F}(ch^2) = c^2(N_{K_1/F}(h))^2$.
    Let $\gamma'=bh$. Then $N_{K/F}(\gamma') =
    c^2(N_{K_1/F}(h))^2$ so that $N_{K/F}(\gamma/\gamma')=1$. By
    Hilbert~90 there exists $\alpha\in K^\times$ with
    $\alpha^{\sigma-1}=\gamma/\gamma'$. Then
    \begin{align*}
        [N_{K/F}(\alpha)] &=[\alpha]^{(\sigma-1)^3}
        =[\gamma/\gamma']^{(\sigma-1)^2}
        =[N_{K/K_1}(\gamma\gamma')] \\ &= [ch^2b^2h^2] = [c]
        =[\gamma]^{(\sigma-1)^2}.
    \end{align*}
    Because $M_\gamma^G=\langle [\gamma]^{(\sigma-1)^2} \rangle$
    our statement follows.
\end{proof}

In what follows, let $l_H(\gamma)$ denote the dimension over $\Fp$
of the cyclic $\Fp[H]$-submodule of $J$ generated by $[\gamma]$.

\begin{lemma}[Second Fixed Elements are Norms Lemma]\label{le:fixednew}

    (a) Suppose $p>2$ and $n\geq 1$. Let $\gamma\in K^\times$ with
    $[\gamma]\in J\setminus [K_{n-1}^\times]$, and let
    $H=\Gal(K/K_{n-1})$. Assume that one of the following holds:
    \begin{itemize}
        \item $\xi_p\notin F$
        \item $\xi_p\in F$ and $l_H(\gamma)\geq 3$
        \item $\xi_p\in F$, $l_H(\gamma)=2$, and
        $[N_{K/F}(\gamma)]_F = [1]_F$.
    \end{itemize}
    Then
    \begin{equation*}
        [\gamma]^{(\sigma-1)^{l(\gamma)-1}} \in [N_{K/F}(K^\times)].
    \end{equation*}

    (b) Suppose $p=2$ and $n\ge 2$.  Let $\gamma\in K^\times$ and
    $H=\Gal(K/K_{n-2})$.  Assume that one of the following holds:
    \begin{itemize}
        \item $l_H(\gamma)=4$
        \item $l_H(\gamma)=3$ and $[N_{K/F}(\gamma)]_F = [1]_F$.
    \end{itemize}
    Then
    \begin{equation*}
        [\gamma]^{(\sigma-1)^{l(\gamma)-1}} \in
        [N_{K/F}(K^\times)].
    \end{equation*}
\end{lemma}

\begin{proof}
    (a). Since part (a) is true for $n=1$ by the First Fixed
    Elements are Norms Lemma~\eqref{le:fixed1}, let us assume that
    $n>1$.  The Fixed Submodule Lemma~\eqref{le:fixed} tells us
    that $l_H(\gamma) \geq 2$, since $[\gamma]\notin
    [K_{n-1}^\times]$.

    Now if $l_H(\gamma)=2$, we claim that $\gamma$ is not
    exceptional for $K/K_{n-1}$, as follows.  Since
    $l_H(\gamma)=2<p$, the Norm Lemma~\eqref{le:norm} tells us that
    $[N_{K/K_{n-1}}(\gamma)]_{K_{n-1}}\in \langle
    [a_{n-1}]_{K_{n-1}}\rangle$.  If $\gamma$ is exceptional for
    $K/K_{n-1}$, then $[N_{K/K_{n-1}}(\gamma)]_{K_{n-1}} \neq
    [1]_{K_{n-1}}$.  By the Norm Lemma~\eqref{le:norm} again,
    $[N_{K/F}(\gamma)]_F \neq [1]_F$, contradicting our hypothesis.
    Hence if $l_H(\gamma)=2$ then $\gamma$ is not exceptional for
    $K/K_{n-1}$, as required.

    Let
    \begin{equation*}
        [\beta] = [\gamma]^{(\sigma^{p^{n-1}}-1)^{l_H(\gamma)-1}}
        = [\gamma]^{(\sigma-1)^{p^{n-1}(l_H(\gamma)-1)}}.
    \end{equation*}
    We invoke the First Fixed Elements are Norms
    Lemma~\eqref{le:fixed1} and deduce that there exists
    $[\alpha]\in J$ such that $[\beta] = [N_{K/K_{n-1}}(\alpha)]$.
    Then
    \begin{equation*}
        [\beta] = [\alpha]^{(\sigma-1)^{p^{n-1}(p-1)}}
    \end{equation*}
    since $l_H(\alpha)=p$. Set $s=l(\beta)$. Then
    \begin{align*}
        [\alpha]^{(\sigma-1)^{p^n-p^{n-1}}(\sigma-1)^{s-1}} &=
        [\beta]^{(\sigma-1)^{s-1}} \\ &=
        [\gamma]^{(\sigma-1)^{p^{n-1}(l_H (\gamma)-1)+s-1}},
    \end{align*}
    and this element is in $J^G$.

    Set $[\lambda] := [\alpha]^{(\sigma-1)^{s}}$. Then we have
    \begin{equation*}
        [\lambda]^{(\sigma-1)^{p^n-p^{n-1}-1}}=
        [\alpha]^{(\sigma-1)^{p^n-p^{n-1}+s-1}} =
        [\beta]^{(\sigma-1)^{s-1}} \neq [1].
    \end{equation*}
    Hence $l(\lambda)=p^n-p^{n-1}$.

    Now we claim that $l_H(\lambda)=p-1$. First, since
    \begin{equation*}
        [\lambda]^{(\sigma^{p^{n-1}}-1)^{p-1}} =
        [\lambda]^{(\sigma-1)^{p^n-p^{n-1}}} = [1]
    \end{equation*}
    we see that $l_H(\lambda)\leq p-1$. But since
    \begin{equation*}
        [\lambda]^{(\sigma-1)^{p^{n-1}(p-2)}} =
        [\lambda]^{(\sigma-1)^{p^n-p^{n-1}-p^{n-1}}}
    \end{equation*}
    and $p^{n-1}\ge 1$, we see that
    \begin{equation*}
        [\lambda]^{(\sigma-1)^{p^{n-1}(p-2)}}\neq [1].
    \end{equation*}
    (In fact, since we assume that $n>1$ we have $p^{n-1}>1$, but we
    do not need the strict inequality.) Therefore indeed
    $l_H(\lambda)=p-1\geq 2$, since we assume that $p\geq 3$.
    Observe that since $[\beta]\neq [1]$ we have $s=l(\beta)>0$.
    Thus $[\lambda]$ is in the image of $\sigma-1$ and hence
    $[N_{K/F}(\lambda)]_F=[1]_F$. Since $l_H(\lambda)=p-1<p$, we
    obtain
    \begin{equation*}
        [N_{K/K_{n-1}}(\lambda)]_{K_{n-1}} \in \langle
        [a_{n-1}]_{K_{n-1}} \rangle.
    \end{equation*}

    By the Norm Lemma~\eqref{le:norm}, we deduce that
    \begin{equation*}
        [N_{K/K_{n-1}}(\lambda)]_{K_{n-1}} = [1]_{K_{n-1}}.
    \end{equation*}
    Hence $\lambda$ is unexceptional for $K/K_{n-1}$. Thus we can
    use the First Fixed Elements are Norms Lemma~\eqref{le:fixed1}
    for $\lambda$. We see that there exists $\chi\in K^\times$
    such that
    \begin{equation*}
        [\lambda]^{(\sigma^{p^{n-1}}-1)^{l_H(\lambda)-1}} =
        [\chi]^{(\sigma-1)^{p^n-p^{n-1}}}
    \end{equation*}
    or equivalently
    \begin{equation*}
        [\lambda]^{(\sigma-1)^{p^n-2p^{n-1}}}=
        [\chi]^{(\sigma-1)^{p^n-p^{n-1}}}.
    \end{equation*}
    This means in particular that
    \begin{equation*}
        l(\chi)=l(\lambda)+p^{n-1}=p^n.
    \end{equation*}
    Putting our calculations together, we obtain
    \begin{align*}
        [N_{K/F}(\chi)] &= [\chi]^{(\sigma-1)^{p^n-1}}
        =[\lambda]^{(\sigma-1)^{p^n-p^{n-1}-1}} \\
        &=[\alpha]^{(\sigma-1)^{p^n-p^{n-1}+s-1}}
        =[\gamma]^{(\sigma-1)^{p^{n-1}(l_H(\gamma)-1)+s-1}} \\
        &=[\gamma]^{(\sigma-1)^{l(\gamma)-1}}
    \end{align*}
    as required.

    (b). If $l_H(\gamma)=3$ then we claim that
    $[N_{K/K_{n-2}}(\gamma)]_{K_{n-2}} = [1]_{K_{n-2}}$, as
    follows. Since $l_H(\gamma)<4$, we have from the Norm
    Lemma~\eqref{le:norm} that $[N_{K/K_{n-2}} (\gamma)]_{K_{n-2}}
    \in \langle [a_{n-2}]_{K_{n-2}} \rangle$.  If $[N_{K/K_{n-2}}
    (\gamma)]_{K_{n-2}} = [a_{n-2}]_{K_{n-2}}^s$ for $s\not\equiv
    0 \bmod 2$, then we obtain from the Norm Lemma~\eqref{le:norm}
    that $[N_{K/F}(\gamma)]_F = [a]_F^s \neq [1]$, contradicting our
    hypothesis.  Therefore if $l_H(\gamma)=3$ then
    $[N_{K/K_{n-2}}(\gamma)]_{K_{n-2}} = [1]_{K_{n-2}}$, as
    required.

    We may then invoke the Fixed Elements of Length 3 Submodules
    are Norms Lemma~\eqref{le:fixed2} and deduce that there exists
    $\alpha \in K^\times$ such that
    \begin{equation*}
        [\alpha]^{(\sigma^{2^{n-2}}-1)^3}=
        [\gamma]^{(\sigma^{2^{n-2}}-1)^{l_H(\gamma)-1}}.
    \end{equation*}

    If instead $l_H(\gamma)=4$, then by setting $\alpha=\gamma$ we
    see that $\alpha$ as above exists as well.

    In either case, then, we obtain the equation with $\alpha$
    above. Hence
    \begin{equation*}
        [\alpha]^{(\sigma-1)^{2^n-2^{n-2}}} =
        [\gamma]^{(\sigma-1)^{2^{n-2}(l_H(\gamma)-1)}}\neq [1].
    \end{equation*}
    Set $s:=l(\gamma)-2^{n-2}(l_H(\gamma)-1)>0$. Then we have
    \begin{equation*}
        [\alpha]^{(\sigma-1)^{2^n-2^{n-2}+s-1}} =
        [\gamma]^{(\sigma-1)^{2^{n-2}(l_H(\gamma)-1)+s-1}}\neq
        [1].
    \end{equation*}
    Furthermore, this element belongs to $J^G$. Set $[\lambda]:=
    [\alpha]^{(\sigma-1)^s}$. Then
    \begin{equation*}
        [\lambda]^{(\sigma-1)^{2^n-2^{n-2}-1}} =
        [\alpha]^{(\sigma-1)^{2^n-2^{n-2}+s-1}},
    \end{equation*}
    whence $l(\lambda)=2^n-2^{n-2}$.

    Now consider $l_H(\lambda)$. On the one hand,
    \begin{equation*}
        [\lambda]^{(\sigma^{2^{n-2}}-1)^3} =
        [\lambda]^{(\sigma-1)^{2^n-2^{n-2}}} = [1],
    \end{equation*}
    and on the other hand
    \begin{equation*}
        [\lambda]^{(\sigma^{2^{n-2}}-1)^2}\neq [1].
    \end{equation*}
    We deduce that $l_H(\lambda)=3$. Observe that since
    $[\lambda]$ is in the image of $\sigma-1$ we have
    $[N_{K/F}(\lambda)]_F=[1]_F$. Since $l_H(\lambda)=3$, we see
    that
    \begin{equation*}
        [N_{K/K_{n-2}}(\lambda)]_{K_{n-2}} \in \langle
        [a_{n-2}]_{K_{n-2}} \rangle.
    \end{equation*}

    By the Norm Lemma~\eqref{le:norm}, we deduce that
    \begin{equation*}
        [N_{K/K_{n-2}}(\lambda)]_{K_{n-2}}=[1]_{K_{n-2}}.
    \end{equation*}
    By the Fixed Elements of Length 3 Submodules are Norms
    Lemma~\eqref{le:fixed2}, there exists $\chi\in K^\times$ with
    \begin{equation*}
        [\chi]^{(\sigma^{2^{n-2}}-1)^3} =
        [\lambda]^{(\sigma^{2^{n-2}}-1)^2}.
    \end{equation*}
    Equivalently
    \begin{equation*}
        [\chi]^{(\sigma-1)^{2^n-2^{n-2}}} =
        [\lambda]^{(\sigma-1)^{2^{n-1}}},
    \end{equation*}
    and therefore $l(\chi)=l(\lambda)+2^{n-2}=2^n$.

    Summarizing, we have obtained
    \begin{align*}
        [N_{K/F}\chi] &= [\chi]^{(\sigma-1)^{2^n-1}}
        =[\lambda]^ {(\sigma-1)^{2^n-2^{n-2}-1}} \\
        &=[\alpha]^{(\sigma-1)^{2^n-2^{n-2} +s-1}}
        = [\gamma]^{(\sigma-1)^{l (\gamma)-1}}
    \end{align*}
    as required.
\end{proof}

\section{Base Cases}\label{se:base}

\begin{proposition}\label{pr:noxn1}
    Theorem~\ref{th:nox} holds for $n=1$.
\end{proposition}

\begin{proof}
    Let $\Ic$ be an $\Fp$-basis for $[N_{K/F}(K^\times)]$.  For
    each $[x]\in \Ic$, we construct a free $\Fp[G]$-module $M(x)$,
    as follows. Choose a representative $x\in F^\times$ for $[x]$
    such that $x\in N_{K/F}(K^\times)$.  Choose $\gamma\in
    K^\times$ such that $x=N_{K/F}(\gamma)$.  Finally let $M(x)$
    be the $\Fp[G]$-submodule of $J$ generated by $[\gamma]$.
    Since $[N_{K/F}(\gamma)]= [\gamma]^{(\sigma-1)^{p-1}}=[x]\neq
    [1]$, $\dim_{\Fp} M(x) = p$ and hence $M(x)$ is free.  By the
    Exclusion Lemma~\eqref{le:excl}, the set of modules $M(x)$,
    $[x]\in \Ic$, is independent.

    Let $Y_1=\oplus_{\Ic} M(x)$.  Then $Y_1$ is a free
    $\Fp[G]$-module with $Y_1^G=[N_{K/F}(K^\times)]$.  Let $Y_0$
    be any complement in $[F^\times]$ of $Y_1^G$.  Clearly $Y_0$
    is a trivial $\Fp[G]$-module.  Since $Y_0^G\cap Y_1^G =
    \{0\}$, $Y_0+Y_1 = Y_0\oplus Y_1$ by the Exclusion
    Lemma~\eqref{le:excl}.  Moreover, $(Y_0+Y_1)^G = [F^\times]$.

    Now set $\tilde J := Y_0+Y_1$.  Then, applying the Inclusion
    Lemma~\eqref{le:incl} with $U=\tilde J$, $V=J$, and $U+V=J$,
    we will deduce that $\tilde J=J$.  Observe first that $(U+V)^G
    = J^G$ which, by the Fixed Submodule Lemma~\eqref{le:fixed},
    is $[F^\times]$.  Since $\tilde J^G = [F^\times]$, we obtain
    $(U+V)^G \subset U$.

    Let $[\gamma] \in J \setminus J^G$.  Then $l(\gamma)\ge 2$. If
    $p=2$ then $[c] = [\gamma]^{(\sigma-1)^{l(\gamma)-1}} =
    [\gamma]^{(\sigma-1)} = N[\gamma]$.  Otherwise, by the First
    Fixed Elements are Norms Lemma~\eqref{le:fixed1}, we obtain
    $[c] = [\gamma]^{ (\sigma-1)^{ l(\gamma)-1}} =
    [N_{K/F}(\alpha)]$ for some $\alpha\in K^\times$.  In any
    case, $[c]\in [N_{K/F}(K^\times)]$. Equivalently, switching
    for the moment to additive notation for convenience, $[c] =
    \sum_{\Ic} c_x[x]$ with almost all $c_x=0$. Now for each
    $[x]$, $M(x)=M_{\omega(x)}$ for some $\omega(x)\in K^\times$
    with $N([\omega(x)])=[x]$.  Hence $[c]=N(\sum
    c_x[\omega(x)])\in Y_1 \subset \tilde J$.  We have shown that
    for every $[\gamma]\in J\setminus J^G$,
    $[\gamma]^{(\sigma-1)^{ l(\gamma)-1}} = N([\alpha])$ for
    $[\alpha]\in Y_1\subset \tilde J$.  Hence we have satisfied
    the hypotheses of the Inclusion Lemma~\eqref{le:incl}, and
    $J\subset \tilde J$, as required.
\end{proof}

\begin{proposition}\label{pr:xn1}
    Theorem~\ref{th:x} holds for $n=1$.  In this case $m\in
    \{-\infty, 0\}$, and $p^m+1$ is the minimal length $l(z)$ for
    all $z\in K^\times$ with $[N_{K/F}(z)]_F\neq [1]_F$.
\end{proposition}

\begin{proof}
    Let $X$ be the cyclic submodule of $J$ generated by the given
    exceptional element $[\delta]$.  Since $\delta=\root{p}\of{a}$
    satisfies $[N_{K/F}(\delta)]_F=[a]_F$ and $[\delta]^{(\sigma-1)}
    = [\xi_p] \in J^G$, we have $m<1$.

    For the case in which $p=2$ and $-1\in N_{K/F}(K^\times)$, let
    $\gamma$ satisfy $N_{K/F}(\gamma)=-1$.  Then set
    $\gamma'=\sqrt{a}\gamma$. We have $N_{K/F}(\gamma')=a$ and
    $[\gamma']^{(\sigma-1)} = [\gamma']^{(1+\sigma)} =
    [N_{K/F}(\gamma')] = [a] = [1] \in [K_{-\infty}^\times]$.
    Hence in this case $\gamma'$ is exceptional and $m=-\infty$.
    By the definition, then, for any exceptional $\delta$ in this
    case, we have $[\delta]^{ (\sigma-1)} = [1]$.

    In any case, by the Exact Sequence Lemma~\eqref{le:exact}(i), we
    have $[\delta]\not\in [F^\times]$. If $m=-\infty$, then $X$ is
    of dimension 1 and hence $X\cap [F^\times] = \{0\}$.
    Furthermore, since $l(z)=0$ implies $[z]=[1]$ and
    $[N_{K/F}(z)]_F = [1]_F$, we see that $p^{m}+1$ is the minimal
    $l(z)$ among $z\in K^\times$ with $[N_{K/F}(z)]_F\neq [1]_F$.
    If, on the other hand, $m=0$, then $X$ is of dimension 2 and by
    the Fixed Submodule Lemma~\eqref{le:fixed}, we have $X^G =
    X^{(\sigma-1)} = X \cap [F^\times]$.  Furthermore, $p^m+1=2$ is
    the minimal $l(z)$ among $z\in K^\times$ with
    $[N_{K/F}(z)]_F\neq [1]_F$, since if $l(z)<2$ for such a $z$
    then $[z]^{\sigma-1}=[1]\in [K_{-\infty}^\times]$, a
    contradiction of $m=0$.

    We proceed to construct $Y_1$. Let $\Ic$ be an $\Fp$-basis for
    $[N_{K/F}(K^\times)]$.  For each $[x]\in \Ic$, we construct a
    free $\Fp[G]$-module $M(x)$, as follows. Choose a
    representative $x\in F^\times$ for $[x]$ such that $x\in
    N_{K/F}(K^\times)$. Choose $\gamma\in K^\times$ such that
    $x=N_{K/F}(\gamma)$.  Finally let $M(x)=M_\gamma$, the
    $\Fp[G]$-submodule of $J$ generated by $[\gamma]$. Since
    $[N_{K/F}(\gamma)]= [\gamma]^{(\sigma-1)^{p-1}}=[x]\neq [1]$,
    $\dim_{\Fp} M(x) = p$ and hence $M(x)$ is free. By the
    Exclusion Lemma~\eqref{le:excl}, the set of modules $M(x)$,
    $[x]\in \Ic$, is independent. Let $Y_1=\oplus_{\Ic} M(x)$.
    Then $Y_1$ is a free $\Fp[G]$-module with
    $Y_1^G=[N_{K/F}(K^\times)]$.

    Now $X^G\cap Y_1^G=\{0\}$, as follows.  Suppose not.  Then since
    $X\cap [F^\times] = \{0\}$ in the case $m=-\infty$, we must have
    $m = 0$.  In particular, by the considerations above for $p=2$,
    we see that $p>2$.  Let $X^G\cap Y_1^G \neq \{0\}$. Then since
    $X^G$ is an $\Fp$-vector space generated by
    $[\delta]^{\sigma-1}$ we see that $[\delta]^{\sigma-1}\in
    X^G\cap Y_1^G$.  Because $Y_1$ is a free $\Fp[G]$-module, there
    exists $[\alpha]\in Y_1$ such that $N[\alpha] =
    [\delta]^{\sigma-1}$. Consider $\delta' =
    \delta/(\alpha)^{(\sigma-1)^{p-2}}$.  Then $[N_{K/F}(\delta')]_F
    \neq [1]_F$ and $[\delta']^{\sigma-1} = [1] \in
    [K_{-\infty}^\times]$, so that $m=-\infty$, a contradiction.
    Because $X^G\cap Y_1^G = \{0\}$, by the Exclusion
    Lemma~\eqref{le:excl} we have $X+Y_1 = X\oplus Y_1$.

    Now let $Y_0$ be any complement in $[F^\times]$ of the
    $\Fp$-submodule of $J$ generated by $X\cap [F^\times]$ and
    $Y_1^G$. Clearly $Y_0$ is a trivial $\Fp[G]$-module. Since
    $Y_0^G\cap (X+Y_1)^G = \{0\}$, we obtain $X + Y_0 + Y_1 = X
    \oplus Y_0 \oplus Y_1$ from the Exclusion
    Lemma~\eqref{le:excl}.

    If $m=-\infty$ then observe that $[F^\times] = Y_0^G + Y_1^G$,
    and if $m=0$ then since $X^G=X^{(\sigma-1)}$, we have
    $[F^\times] = X^{(\sigma-1)} + Y_0^G + Y_1^G$.

    Now set $\tilde J=X+Y_0+Y_1$.  We adapt the proof of the
    Inclusion Lemma~\eqref{le:incl} to show that $J\subset \tilde
    J$ and hence $J=\tilde J$, by induction on the socle series
    $J_i$ of $J$.

    We first show that if $[\beta]\in J_1=J^G$ then $[\beta]\in
    \tilde J$. If $[N_{K/F}\beta]_F = [1]_F$, then the Proper
    Subfield Lemma~\eqref{le:prop} gives $[\beta] \in [F^\times]$.
    Since $Y_0$ is a complement in $[F^\times]$ of the submodule
    generated by $X \cap [F^\times]$ and $Y_1^G$, $[\beta]\in \tilde
    J$.

    Otherwise $[N_{K/F}\beta]_F \neq [1]_F$. Since $l(\beta)=1$ we
    must have $m=-\infty$ and $[\delta]\in J_1$.  By the Exact
    Sequence Lemma~\eqref{le:exact}, both $[N_{K/F}(\beta)]_F$ and
    $[N_{K/F}(\delta)]_F$ lie in $\langle [a]_F \rangle$, and by the
    definition of exceptionality, both are generators of $\langle
    [a]_F \rangle$.  Hence $[N_{K/F}(\beta)]_F =
    [N_{K/F}(\delta)]_F^s$ for some $s\in \Z$, and we set
    $\beta'=\beta/\delta^s$.  Then $[\beta']\in J^G$ and
    $[N_{K/F}(\beta')]_F = [1]_F$.  By the Exact Sequence
    Lemma~\eqref{le:exact}(ii), we see that $[\beta']\in
    [F^\times]$. As in the preceding paragraph, this gives
    $[\beta']\in \tilde J$. Then, since $[\delta]\in \tilde J$ as
    well, we obtain $[\beta]\in \tilde J$.  Hence $J_1\subset \tilde
    J$.

    For the inductive step, assume that $J_i\subset \tilde J$ for
    all $1\le i<t\le p$, and let $[\gamma] \in J_t \setminus
    J_{t-1}$.

    We first claim that in the particular case of $t=2$, without
    loss of generality we may assume that $\gamma$ is unexceptional,
    as follows.  Assume that $\gamma$ is exceptional and
    $l(\gamma)=2$.  Then $m\neq -\infty$ since otherwise $l(\gamma)$
    would be $1$.  Since $n=1$ we have $m=0$. We established
    earlier, however, that if $p=2$ and $n=1$ then $m=-\infty$.
    Hence $p>2$. Now since $m=0$, we have $[\delta]^{(\sigma-1)} \in
    [F^\times]$, $l(\delta)\le 2$, and by the definition of
    exceptionality $l(\delta) \neq 1$, since otherwise $m=-\infty$.
    Hence $l(\delta)=l(\gamma)=2$.

    Since $p>2$, we deduce from the Norm Lemma~\eqref{le:norm}
    that both of $[N_{K/F}(\gamma)]_F$ and $[N_{K/F}(\delta)]_F$
    lie in $\langle [a]_F \rangle$, and by the definition of
    exceptionality, both are generators.  Hence
    $[N_{K/F}(\gamma)]_F = [N_{K/F}(\delta)]_F^s$ for some $s\in
    \Z$, and we set $\gamma'=\gamma/\delta^s$.  Then
    $l(\gamma')\le 2$ and $[N_{K/F}(\gamma')]_F = [1]_F$, so that
    $\gamma'$ is unexceptional.  Since $[\delta]\in X\subset
    \tilde J$, to show that $[\gamma]\in \tilde J$ it is enough to
    show that $[\gamma']\in J$.  We may therefore assume that
    $\gamma$ is unexceptional if $t=2$.

    Now if $p=2$ then $l(\gamma)=2$ and
    \begin{equation*}
        [c] = [\gamma]^{(\sigma-1)^{l(\gamma)-1}} =
        [\gamma]^{(\sigma+1)} = N[\gamma] = [N_{K/F}(\gamma)],
    \end{equation*}
    and we set $\alpha=\gamma$.  Otherwise, $p>2$ and by the First
    Fixed Elements are Norms Lemma~\eqref{le:fixed1}, we have $[c]
    = [\gamma]^{ (\sigma-1)^{ l(\gamma)-1}} = [N_{K/F}(\alpha)]$
    for some $\alpha\in K^\times$.

    In either case, $[c]\in [N_{K/F}(K^\times)]$. Equivalently,
    switching for the moment to additive notation for convenience,
    $[c] = \sum_{\Ic} c_x[x]$ with almost all $c_x=0$. Now for
    each $[x]$, $M(x) = M_{\omega(x)}$ for some $\omega(x)\in
    K^\times$ with $N([\omega(x)])=[x]$.  Hence $[c] = N(\sum
    c_x[\omega(x)])\in Y_1 \subset \tilde J$.  Switching back to
    multiplicative notation, $[c] = [\alpha]^{(\sigma-1)^{p-1}}$
    for some $[\alpha]\in Y_1$. Let $[\gamma'] =
    [\alpha]^{(\sigma-1)^{p-t}} \in \tilde J$. Since
    $[\gamma/\gamma']^{(\sigma-1)^{t-1}} = [1]$, we find
    $l(\gamma/\gamma')<t$.  By induction, $[\gamma/\gamma'] \in
    \tilde J$, and hence $[\gamma]\in \tilde J$ as well.  By
    induction on the socle series, then, $J\subset \tilde J$.

    Finally we verify part (3) of Theorem~\ref{th:x} in our case
    $n=1$.  In this case $i=0$ and $m\le 0$.  If $m=-\infty$ we have
    $X^{(\sigma-1)}= \{[1]\}$ and
    \begin{equation*}
        X^{(\sigma-1)}\oplus Y^{H_0} = Y^G = [F^\times] =
        [K_0^\times]
    \end{equation*}
    by the construction of $Y=Y_0\oplus Y_1$.

    Suppose instead that $m=0$.  Then $X\cap [F^\times] =
    X^{(\sigma-1)}$ and again from our construction of $Y_0$ and
    $Y_1$ we obtain
    \begin{equation*}
        X^{(\sigma-1)} \oplus Y_0 \oplus Y_1^G = X^{(\sigma-1)}
        \oplus Y^G = [F^\times] = [K_0^\times]
    \end{equation*}
    as required.
\end{proof}

\begin{proposition}\label{pr:xp2n2}
    Theorem~\ref{th:x} holds in the case $p=2$, $n=2$.  In this case
    $m\in \{-\infty, 0, 1\}$, and $2^m+1$ is the minimal length
    $l(z)$ for all $z\in K^\times$ with $[N_{K/F}(z)]_F\neq [1]_F$.
\end{proposition}

\begin{proof}
    Let $X$ be the cyclic submodule of $J$ generated by the given
    exceptional element $[\delta]$.  Consider $\theta=\sqrt{a_1}$.
    Then $[N_{K/F}(\theta)]_F=[a]_F$.  Because $K/F$ is Galois we
    have $a_1^\sigma=a_1k^2$ for some $k\in K_1^\times$.
    Therefore $[\theta]^{\sigma-1} = [\pm k]\in K_1^\times$.
    Hence $m<2$.  Now let $\delta$ be any exceptional element in
    $K^\times$.  Because $[N_{K/F}(\delta)]_F\neq [1]_F$ we see
    that $[\delta]\not\in [F^\times]$.

    If $m=-\infty$, then $X$ is of dimension 1 and therefore $X\cap
    [F^\times] = \{0\}$. If $m=0$, then $X$ is of dimension 2, and
    by the Exact Sequence Lemma~\eqref{le:exact}(i,ii), observing
    that $N_{K/F}(\delta^{\sigma-1})=1$, we obtain $X^G =
    X^{(\sigma-1)} = X \cap [F^\times]$.  Finally assume that $m=1$.
    Observe that then $l(\delta^{\sigma-1}) \neq 1$. Indeed
    otherwise $N_{K/F}(\delta^{\sigma-1}) = 1$ and the Exact
    Sequence Lemma~\eqref{le:exact}(ii) implies that
    $[\delta]^{\sigma-1} \in [F^\times]$, which contradicts our
    assumption that $m=1$. Hence $l(\delta)\geq 3$. However since
    $[\delta]^{\sigma -1} \in [K_1^\times]$ and $(\sigma-1)^2\equiv
    \sigma^2-1$, we have $l(\delta^{\sigma-1}) \leq 2$, and
    therefore $l(\delta) \leq 3$. Consequently $l(\delta)=3$. Since
    $[N_{K/F}(\delta)]_F \neq [1]_F$ we see that $[\delta] \notin
    [K_1^\times]$. Therefore $X^{(\sigma-1)}=X\cap [K_1^\times]$.

    We claim that $2^m+1$ is the minimal $l(z)$ among $z\in
    K^\times$ with $[N_{K/F}(z)]_F\neq [1]_F$, as follows.  If
    $l(z)=0$ then certainly $[z]=[1]$ and $[N_{K/F}(z)]_F = [1]_F$, so
    $l(z)\ge 1$ for any $z$ with $[N_{K/F}(z)]_F\neq [1]_F$.  If
    $m=-\infty$ the claim follows immediately.  If $m=0$, observe
    that if $l(z)<2$ for some $z$ with $[N_{K/F}(z)]_F \neq [1]_F$,
    then $[z]^{\sigma-1}=[1]\in [K_{-\infty}^\times]$, a
    contradiction of $m=0$.  Finally, if $m=1$ then we similarly
    have $l(z)\ge 2$. If further $l(z)=2<2^1+1=3$, then
    $l(z^{\sigma-1})=1$, and the Exact Sequence
    Lemma~\eqref{le:exact}(ii) implies that $[z]^{\sigma-1} \in
    [F^\times]$, a contradiction of $m=1$.

    We proceed to construct $Y_2$. Let $\Ic_2$ be an $\Ft$-basis
    for $[N_{K/F}(K^\times)]$.  For each $[x]\in \Ic_2$, we
    construct a free $\Ft[G]$-module $M(x)$, as follows. Choose a
    representative $x\in F^\times$ for $[x]$ such that $x\in
    N_{K/F}(K^\times)$. Choose $\gamma\in K^\times$ such that
    $x=N_{K/F}(\gamma)$. Finally let $M(x)=M_\gamma$. Since
    $[N_{K/F}(\gamma)]= [\gamma]^{(\sigma-1)^{3}} = [x]\neq [1]$,
    $\dim_{\Ft} M(x) = 4$ and hence $M(x)$ is free. By the
    Exclusion Lemma~\eqref{le:excl}, the set of modules $M(x)$,
    $[x]\in \Ic_2$, is independent. Let $Y_2=\oplus_{\Ic_2} M(x)$.
    Then $Y_2$ is a free $\Ft[G]$-module with
    $Y_2^G=[N_{K/F}(K^\times)]$.

    Suppose $X^G\cap Y_2^G\neq \{0\}$. Since $Y_2^G\subset
    [F^\times]$ and $X\cap [F^\times]=\{0\}$ if $m=-\infty$, we are
    in the case $m=0$ or $m=1$, and $X^G=X^{(\sigma-1)^{m+1}} =
    X\cap [F^\times]$. In particular, $l(\delta)=m+2\le 3$. Let
    $f\in F^\times$ satisfy $[1]\neq [f]\in X^G\cap Y_2^G$. Since
    $Y_2$ is free, there exists $[\alpha]\in Y_2$ such that
    $N[\alpha]=[f]$. Let $\delta' =
    \delta/(\alpha)^{(\sigma-1)^{4-l(\delta)}}$.  Then
    $[N_{K/F}(\delta')]_F = [N_{K/F}(\delta)]_F$ since
    $\alpha^{(\sigma-1)^{4-l(\delta)}}$ is in the image of
    $\sigma-1$. Moreover, $l(\delta')<l(\delta)$.  If $m=0$ then
    $l(\delta')\le 1$ and by the definition of exceptionality,
    $m=-\infty$, a contradiction.  If $m=1$ then $l(\delta')\le 2$
    so that $l((\delta')^{\sigma-1})\le 1$ and
    $[(\delta')^{\sigma-1}]\in J^G$. But since
    $(\delta')^{\sigma-1}$ is in the image of $\sigma-1$, we have
    $[N_{K/F}((\delta')^ {\sigma-1})]_F = [1]_F$, and from the Exact
    Sequence Lemma~\eqref{le:exact}(ii) we obtain $[\delta']\in
    [F^\times]$. Then by the definition of exceptionality, $m\le 0$,
    again a contradiction. Thus $X^G\cap Y_2^G=\{0\}$.

    Because $X^G\cap Y_2^G = \{0\}$, by the Exclusion
    Lemma~\eqref{le:excl} we have $X+Y_2 = X\oplus Y_2$.

    We proceed now to construct $Y_1$. Let $\Ic_1$ be an $\Ft$-basis
    for a complement in $[N_{K_1/F} (K_1^\times)]$ of the
    $\Ft$-submodule generated by $[N_{K/F}(K^\times)]$ and $X \cap
    [N_{K_1/F}(K_1^\times)]$. For each $[x]\in
    \Ic_1$, we construct an $\Ft[G]$-module $M(x)$ of dimension 2,
    as follows. Choose a representative $x\in F^\times$ for $[x]$
    such that $x\in N_{K_1/F}(K_1^\times)$. Choose $\gamma\in
    K_1^\times$ such that $x=N_{K_1/F}(\gamma)$.  Finally let
    $M(x)=M_\gamma$. Since $[N_{K_1/F}(\gamma)]=
    [\gamma]^{(\sigma-1)} = [x]\neq [1]$, $\dim_{\Ft} M(x) = 2$.
    The $M(x)$, $[x]\in \Ic_1$, are independent as above. Let
    $Y_1=\oplus_{\Ic_1} M(x)$. Then $Y_1$ is a direct sum of
    $\Ft[G]$-modules of dimension 2, and $Y_1^G$ is the $\Ft$-span
    of $\Ic_1$.  By construction $Y_1^G \cap Y_2^G = \{0\}$ and
    hence by the Exclusion Lemma~\eqref{le:excl}, we have
    $Y_1+Y_2=Y_1\oplus Y_2$.

    Suppose $X^G\cap (Y_1+Y_2)^G \neq \{0\}$. Since $(Y_1+Y_2)^G
    \subset [F^\times]$ and $X\cap [F^\times]=\{0\}$ if
    $m=-\infty$, we are in the case $m=0$ or $m=1$, and $X^G =
    X^{(\sigma-1)^{m+1}} = X\cap [F^\times]$.  Let $X^G=\langle
    [x]\rangle$; then $[1]\neq [x]=[y_1]+[y_2]$ for some $[y_1]\in
    Y_1^G$ and $[y_2]\in Y_2^G$. Since $Y_1^G + Y_2^G \subset
    [N_{K_1/F}(K_1^\times)]$, we deduce $[x]\in
    [N_{K_1/F}(K_1^\times)]$.  We have already established that
    $[y_1]\neq [1]$, since $X^G\cap Y_2^G\neq \{0\}$.  Hence
    $[1]\neq [y_1] = [y_2] + [x]$. But then $Y_1^G$ does not
    consist of a complement of the $\Ft$-submodule generated by
    $Y_2^G=[N_{K/F}(K^\times)]$ and $X\cap
    [N_{K_1/F}(K_1^\times)]$, a contradiction.  Hence we have
    established our equality $X^G\cap (Y_1+Y_2)^G = \{0\}$.

    Because $X^G\cap (Y_1+Y_2)^G=\{0\}$, by the Exclusion
    Lemma~\eqref{le:excl} we have $X+Y_1+Y_2=X\oplus Y_1\oplus
    Y_2$.

    Finally let $Y_0$ be any complement in $[F^\times]$ of the
    $\Ft$-submodule of $J$ generated by $X\cap [F^\times]$,
    $Y_1^G$, and $Y_2^G$. Clearly $Y_0$ is a trivial
    $\Ft[G]$-module. Since $Y_0^G\cap (X+Y_1+Y_2)^G = \{0\}$, we
    see that $X+Y_0+Y_1+Y_2 = X\oplus Y_0\oplus Y_1 \oplus Y_2$ by
    the Exclusion Lemma~\eqref{le:excl}.

    If $m=-\infty$ then observe that $[F^\times] = Y_0^G + Y_1^G +
    Y_2^G$, and otherwise since $X^G=X^{(\sigma-1)^{m+1}}$, we
    have $[F^\times] = X^{(\sigma-1)^{m+1}} + Y_0^G + Y_1^G +
    Y_2^G$.  In order to connect this expression with
    Theorem~\ref{th:x}, part (3), in the case $i=0$, observe
    that $(\sigma -1) (\sigma-1)^{2^m-1} = (\sigma-1)^{2^m} =
    (\sigma-1)^{m+1}$ for $m=0$ or $1$.

    Now let $\tilde J=X+Y_0+Y_1+Y_2$.  We show that $J=\tilde J$
    by showing that an arbitrary element $[\beta]\in J$ lies in
    $\tilde J$, as follows.

    First, if $\beta$ is exceptional, then since $m\le 1$ we have
    $[\beta]^{\sigma-1} \in [K_1^\times]$.  Since $1+\sigma\equiv
    \sigma-1$ on $J$ and $[N_{K_1/F}(\gamma)] =
    [\gamma]^{\sigma+1}$ for $\gamma \in K_1^\times$, we see that
    $l(\beta) \le 3$.  By the Norm Lemma~\eqref{le:norm}, we have
    $[N_{K/F}(\beta)]_F=[a]_F^s$ for some $s \not\equiv 0 \bmod
    p$. Because $p=2$ and $[N_{K/F}(\beta)]_F \neq [1]_F$ we have
    $[N_{K/F}(\beta)]_F = [a]_F$.  Since $\delta$ is exceptional,
    $[N_{K/F}(\delta)]_F=[a]_F$ as well. Then
    $\beta'=\beta/\delta$ satisfies $[N_{K/F}(\beta')]_F =
    [1]_F$ and is therefore unexceptional. Since $[\delta]\in
    X\subset \tilde J$, to show that $[\beta]\in \tilde J$ it
    suffices to show that $[\beta']\in \tilde J$. Therefore we may
    and do assume that $[\beta]$ is unexceptional.

    Observe that the above argument applies not only to elements
    $\beta$ that are exceptional, but in fact to all elements
    $\beta$ such that $[N_{K/F}(\beta)]_F = [a]_F^s$ for some $s
    \in \Z$. Therefore we may assume not only that $\beta$ is
    unexceptional, but also that $[N_{K/F}(\beta)]_F = [1]_F$.

    Suppose that $l(\beta)=1$ and $[N_{K/F}(\beta)]_F = [1]_F$.
    From the Exact Sequence Lemma~\eqref{le:exact}(ii) we see that
    $[\beta] \in [F^\times]$. Since $[F^\times] \subset \tilde J$,
    we obtain $[\beta] \in \tilde J$ as well.

    Now if $l(\beta)=2$, then $[\beta]^{(\sigma^2-1)} =
    [\beta]^{(\sigma-1)^2} = [1]$ and $[\beta]\in J^{H_1}$.
    Moreover, we assume that $[N_{K/F}(\beta)]_F=[1]_F$. By the
    Proper Subfield Lemma~\eqref{le:prop}, we deduce that
    $[\beta]\in [K_1^\times]$. Hence we may assume that the
    representative $\beta$ of $[\beta]$ lies in $K_1^\times$. Then
    $[\beta]^{(\sigma-1)} = [N_{K_1/F}(\beta)] \subset
    [N_{K_1/F}(K_1^\times)]$.

    If $m=1$ then
    \begin{equation*}
        [N_{K_1/F}(K_1^\times)] \subset
        X^{(\sigma-1)^2}+Y_1^G+Y_2^G = X^{(\sigma-1)^2} +
        Y_1^{(\sigma-1)} + Y_2^{(\sigma-1)^3},
    \end{equation*}
    since $Y_1$ is a direct sum of cyclic modules of length 2 and
    $Y_2$ a direct sum of cyclic modules of length 4.
    If $m=0$ then $X \cap [F^\times] = X^{\sigma-1}$ and therefore
    \begin{equation*}
        [N_{K_1/F}(K_1^\times)] \subset X^{(\sigma-1)} +
        Y_1^{(\sigma -1)}+Y_2^{(\sigma-1)^3}.
    \end{equation*}
    If $m=-\infty$ then $X\cap [N_{K_1/F}(K_1^\times)] = \{0\}$
    and
    \begin{equation*}
        [N_{K_1/F}(K_1^\times)] \subset Y_1^G+Y_2^G = Y_1^{(\sigma
        -1)} + Y_2^{(\sigma-1)^3}.
    \end{equation*}
    In any case, $[\beta]^{(\sigma-1)}$ lies in $\tilde
    J^{\sigma-1}$ and hence there exists $\alpha\in \tilde J$ such
    that $[\alpha]^{(\sigma-1)} = [\beta]^{(\sigma-1)}$.  But then
    $[\alpha/\beta]\in J^G$, which we have already established
    lies in $\tilde J$. Hence $[\beta]\in \tilde J$.

    Now suppose that $l(\beta)\ge 3$ and $[N_{K/F}(\beta)]_F =
    [1]_F$. By the Fixed Elements of Length 3 Submodules are Norms
    Lemma~\eqref{le:fixed2}, we have $[c] = [\beta]^{ (\sigma-1)^{
    l(\beta)-1}} = N[\alpha]$ for some $\alpha\in K^\times$.
    Equivalently, switching for the moment to additive notation for
    convenience, $[c] = \sum_{\Ic_2} c_x[x]$ with almost all
    $c_x=0$. As in the proof of the previous proposition, we obtain
    $[c]=N(\sum c_x[\omega(x)])\in Y_2 \subset \tilde J$. Let
    $[\beta'] = [\beta]-(\sigma-1)^{4-l(\beta)}(\sum
    c_x[\omega(x)])$. Then $l(\beta')<l(\beta)$ and we proceed by
    induction.

    Hence $J=\tilde J$.

    Now we consider the location of $[K_1^\times]$ in $J$.  Since
    $K_1$ is the fixed field in $K$ of $H_1$, we have
    \begin{equation*}
        [K_1^\times]\subset J^{H_1} = X^{H_1} \oplus Y_0 \oplus Y_1
        \oplus Y_2^{H_1}.
    \end{equation*}
    By our construction of $Y_0$ and $Y_1$ we see that $Y_0 \oplus
    Y_1 \subset [K_1^\times]$. By the Submodule-Subfield
    Lemma~\eqref{le:sub} we see that $Y_2^{H_1}= Y_2 \cap
    [K_1^\times]$. Also because $m\leq 1$ we see from the
    definition of $m$ that $X^{(\sigma-1)} \subset [K_1^\times]$.
    Hence $X^{(\sigma-1)} + Y_0 + Y_1 + Y_2^{H_1} \subset
    [K_1^\times]$.  It remains to show that this inclusion is an
    equality.

    We showed after the definition of exceptional element
    that $[\delta] \notin [K_{n-1}^\times] = [K_1^\times]$. Therefore
    $X \cap [K_1^\times] = X^{(\sigma-1)}$, and
    we have
    \begin{equation*}
        X^{(\sigma-1)}\oplus Y_0 \oplus Y_1 \oplus Y_2^{H_1}
        \subset [K_1^\times] \subset X^{H_1} \oplus Y_0 \oplus Y_1
        \oplus Y_2^{H_1}.
    \end{equation*}
    Hence each $[k] \in [K_1^\times]$ can be written as
    \begin{equation*}
        [k]=[x]+[y],\text{ where }[x]\in X^{H_1}
        \text{ and }[y]\in Y_0 \oplus Y_1 \oplus Y_2^{H_1}.
    \end{equation*}
    Thus
    \begin{equation*}
        [x]=[k]+[y] \in [K_1^\times] + Y_0 + Y_1 + Y_2^{H_1} \subset
        X \cap [K_1^\times].
    \end{equation*}
    Therefore we see that $X^{(\sigma-1)}\oplus Y_0 \oplus Y_1
    \oplus Y_2^{H_1} = [K_1^\times]$. Observe that if $m=-\infty$
    then $X^{(\sigma-1)} = \{0\}$. Since $m\leq 1$ we see that our
    decomposition of $[K_1^\times]$ is in agreement with
    Theorem~\ref{th:x}, part (3).
\end{proof}

\section{Free Submodules and Proof of Theorem~1}
\label{se:nox}

For the following proposition, assume that Theorems~\ref{th:nox} and
\ref{th:x} hold for all extensions of degree $p^s$, $1\le s< n$, and
if $p=2$, then $n>2$.

\begin{proposition}\label{pr:y}
    There exists a submodule
    \begin{equation*}
        \hat Y = \hat Y_n \oplus \hat Y_{n-1} \oplus \dots
        \oplus \hat Y_0
    \end{equation*}
    of $J$ such that
    \begin{enumerate}
        \item $\hat Y_i$ is a direct sum of cyclic
        $\Fp[G]$-modules of dimension $p^i$;
        \item $[K_i^\times]=\hat Y^{H_i}$ for $0 \le i < n$;
        \item $\hat Y_n^G = [N_{K/F}(K^\times)]$.
     \end{enumerate}
\end{proposition}

\begin{proof}
    Let $\Ic$ be an $\Fp$-base for $[N_{K/F}(K^\times)]$.  As
    usual, for each $[x]\in \Ic$ construct free independent
    $\Fp[G]$-modules $M(x)$, $[x]\in \Ic$, such that $M(x)^G=
    \left<[x]\right>$. Set $\hat Y_n=\oplus_{[x]\in \Ic}
    M(x)$.  Hence $\hat Y_n$ is a direct sum of cyclic
    $\Fp[G]$-modules of dimension $p^{n}$, and $\hat Y_n^G =
    [N_{K/F}(K^\times)]$.

    Assume now that $\xi_p\in F^\times$.  Since $K_{n-1}/F$ embeds
    in a cyclic extension $K$ of degree $p^n$ over $F$,
    $[a_{n-1}]_{K_{n-1}}^{\bar \sigma}=[a_{n-1}]_{K_{n-1}}$ by
    Kummer theory, where $\bar\sigma\in G/H_{n-1}$ is the image of
    $\sigma$ under the natural projection $G\to \bar G :=
    G/H_{n-1}$. (Indeed since $a_{n-1}$ is a $p$th power in $K$,
    so is $a_{n-1}^{\bar \sigma}$; therefore by Kummer theory
    $[a_{n-1}]_{K_{n-1}}^{\bar \sigma} \in \langle
    [a_{n-1}]_{K_{n-1}}\rangle$. However, viewing $\langle
    [a_{n-1}]_{K_{n-1}}\rangle$ as $\Fp$, then $\bar \sigma$ is an
    exponent $p^{n-1}$ action on $\Fp$. Since
    \begin{equation*}
        \Aut(\Fp)\cong \Z/(p-1)\Z,
    \end{equation*}
    this action must be the identity. Hence
    $[a_{n-1}]_{K_{n-1}}^{(\bar \sigma-1)} = [1]_{K_{n-1}}$.)
    Moreover, we have $[N_{K_{n-1}/F}(a_{n-1})]_F = [a]_F$ by
    Proposition~\ref{pr:subgen}.

    Because Theorem~\ref{th:x} holds for $n-1$, we have an
    $\Fp[\bar G]$-module decomposition
    \begin{equation*}
        J(K_{n-1}) = K_{n-1}^\times/K_{n-1}^{\times p} = \langle
        [a_{n-1}]_{K_{n-1}} \rangle \oplus \tilde Y_{n-1}\oplus
        \dots \oplus \tilde Y_0
    \end{equation*}
    into direct sums $\tilde Y_i$ of cyclic $\Fp[\bar G]$-modules
    of dimension $p^i$ and a $\bar G$-invariant submodule $\langle
    [a_{n-1}]_{K_{n-1}} \rangle_{\Fp}$.  Indeed we only have to
    check that $a_{n-1}$ is an exceptional element in $K_{n-1}$.
    This follows since we have shown both
    $[N_{K_{n-1}/F}(a_{n-1})]_F = [a]_F$ and
    $[a_{n-1}]_{K_{n-1}}^{\bar\sigma-1} = [1]_{K_{n-1}}$.

    Moreover, by the Submodule-Subfield Lemma~\eqref{le:sub}
    \begin{equation*}
        \tilde Y_{n-1}^{\bar G} = \tilde Y_{n-1} \cap [N_{K_{n-1}/F}
        (K_{n-1}^\times)]_{K_{n-1}}.
    \end{equation*}
    Because $N_{K_{n-1}/F}$ acts on $J(K_{n-1})$ as
    $(\bar\sigma-1)^{p^{n-1}-1}$ we see that $N_{K_{n-1}/F}$
    annihilates the sum $\tilde Y_{n-2}\oplus \dots \oplus \tilde
    Y_0$.  Also $[N_{K_{n-1}/F}(a_{n-1})]_{K_{n-1}} =
    [1]_{K_{n-1}}$. Therefore
    \begin{equation*}
        [N_{K_{n-1}/F}(K_{n-1}^\times)]_{K_{n-1}} = \tilde Y^{\bar
        G}_{n-1}.
    \end{equation*}

    Assume now that $\xi_p\not\in F^\times$.  Then because
    Theorem~\ref{th:nox} holds for $n-1$, we have an $\Fp[\bar
    G]$-module decomposition
    \begin{equation*}
        J(K_{n-1}) = K_{n-1}^\times/K_{n-1}^{\times p} = \tilde
        Y_{n-1}\oplus \dots \oplus \tilde Y_0
    \end{equation*}
    into direct sums $\tilde Y_i$ of cyclic $\Fp[\bar G]$-modules
    of dimension $p^i$.  As before let $\bar\sigma$ denote the
    image of $\sigma$ under the natural projection $G\to \bar G$.
    Because $N_{K_{n-1}/F}$ acts on $J(K_{n-1})$ as
    $(\bar\sigma-1)^{p^{n-1}-1}$ we see that $N_{K_{n-1}/F}$
    annihilates the sum $\tilde Y_{n-2}\oplus \dots \oplus \tilde
    Y_0$.  Therefore again
    \begin{equation*}
        [N_{K_{n-1}/F}(K_{n-1}^\times)]_{K_{n-1}} = [N_{K_{n-1}/F}
    (\tilde Y_{n-1})]_{K_{n-1}}=\tilde Y^{\bar G}_{n-1}.
    \end{equation*}

    In both cases $\xi_p\in F^\times$, $\xi_p\not\in F^\times$,
    consider $J$ as an $\Fp[H_{n-1}]$-module. Then the Exact
    Sequence Lemma~\eqref{le:exact} gives us that the image of
    each $\tilde Y_0, \dots, \tilde Y_{n-1}$ under the map
    \begin{equation*}
        \epsilon \colon J(K_{n-1}) \to J(K)
    \end{equation*}
    is a direct sum of modules of dimension $p^i$ and that the
    images are all independent. Because the modules $\tilde Y_i$ are
    cyclic as $\Fp[\bar G]$-modules, the images $\epsilon(\tilde
    Y_i)$ are cyclic as $\Fp[G]$-modules. Set $\hat
    Y_i=\epsilon(\tilde Y_i)$ for $i<n-1$. (Recall that we already
    defined $\hat Y_n$ at the beginning of our proof.)

    Set $W := \hat Y_n^{H_{n-1}}$. By the Submodule-Subfield
    Lemma~\eqref{le:sub}
    \begin{equation*}
        W = \hat Y_n^{(\sigma-1)^{p^n-p^{n-1}}} = \hat Y_n \cap
        [K_{n-1}^\times].
    \end{equation*}
    Since $\hat Y_n$ is a direct sum of cyclic $\Fp[G]$-modules of
    dimension $p^n$, $W$ is a direct sum of cyclic modules of
    dimension $p^{n-1}$ and hence is free as an $\Fp[\bar
    G]$-module. Because $W\subset [K_{n-1}^\times]$, we may
    consider the image $P$ of the projection map $\pr \colon W \to
    \epsilon(\tilde Y_{n-1})$ from $W$ to the summand
    $\epsilon(\tilde Y_{n-1})$ in the decomposition
    \begin{equation*}
        [K_{n-1}^\times]  = \epsilon(J(K_{n-1})) = \epsilon(\tilde
        Y_{n-1}) \oplus \hat Y_{n-2} \oplus \cdots \oplus \hat
        Y_0.
    \end{equation*}

    Observe that $W\cong P$ as $\Fp[G]$-modules.  Indeed, since
    $W$ is a free $\Fp[\bar G]$-module, each $[w]\in W\setminus
    \{0\}$ may be written as $[\tilde w]^{(\bar\sigma-1)^s}$ for
    some $0\le s\le p^{n-1}-1$ and $[\tilde w]\in W$ with
    $l(\tilde w)=p^{n-1}$.  We have
    \begin{equation*}
        \pr([\tilde w])^{(\bar\sigma-1)^{p^{n-1}-1}} =
        [\tilde w]^{(\bar\sigma-1)^{p^{n-1}-1}} \neq [1],
    \end{equation*}
    since all other components of $[\tilde w]$ are killed by
    $(\bar\sigma-1)^{p^{n-1}-1}$.  (Since $n\ge 2$, $p^{n-1}-1 \ge
    p^{n-2}$.)  Therefore $\pr([\tilde w])^{(\bar\sigma-1)^{s}} =
    \pr([w]) \neq [1]$.  We conclude that the kernel of the
    projection map is $[1]$, as required.

    Since $M^{\bar G}=M^{(\sigma-1)^{p^{n-1}-1}}$ for free
    $\Fp[\bar G]$-modules, we have further obtained that $W^{\bar
    G} = P^{\bar G}$; equivalently, $W^G=P^G$.  Observe that
    \begin{equation*}
        W^G = W^{\bar G} = W^{(\sigma-1)^{p^{n-1}-1}} \subset
        [N_{K_{n-1}/F}K_{n-1}^\times] =
        \epsilon(\tilde{Y}_{n-1})^G.
    \end{equation*}

    By the Free Complement Lemma~\eqref{le:free}, there exists a
    free $\Fp[\bar G]$-module complement $\hat Y_{n-1}$ in
    $\epsilon(\tilde Y_{n-1})$ of $P$.  Since $W=\hat Y_n\cap
    [K_{n-1}^\times]$, we obtain $\hat Y_n^G = W^{G} = P^{G}$. Now
    the next idea is to use the fact that $\hat Y_n^G=P^G$ to show
    that $\hat Y_n$ and $\hat Y_{n-1}$ are independent and $\hat Y_n
    \oplus \hat Y_{n-1}$ and $\hat Y_{n-2} \oplus \cdots \oplus \hat
    Y_0$ are also independent. Then from the definition of $\hat
    Y_n$ and from our observation above on $\hat Y_i$, $i
    \in\{n-1,\dots,0\}$, it follows immediately that $\hat Y
    =\hat Y_n\oplus\cdots\oplus \hat Y_0 \subset J$ satisfies
    conditions (1) and (3) of our proposition. The last part of our
    proof is then devoted to proving condition (2).

    By the Exclusion Lemma~\eqref{le:excl}, $P^{G}\cap \hat
    Y_{n-1}^{G} = \{0\}$ implies that $\hat Y_{n-1} + \hat Y_n =
    \hat Y_{n-1} \oplus \hat Y_n$. Then, since $P^G + \hat
    Y_{n-1}^G = \epsilon(\tilde Y_{n-1})^G$, we obtain $(\hat
    Y_{n-1} + \hat Y_n)^{G}= \epsilon(\tilde Y_{n-1})^G$. Finally,
    by the Exclusion Lemma~\eqref{le:excl}, $\hat Y_{n-1}+ \hat
    Y_n$ is independent from $\hat Y_{n-2} + \dots + \hat Y_0$.
    Hence we have a submodule
    \begin{equation*}
        \hat Y = \hat Y_n\oplus \hat Y_{n-1} \oplus \dots \oplus
        \hat Y_0 \subset J
    \end{equation*}
    satisfying items (1) and (3).

    We turn next to item (2) and prove that $\hat Y^{H_{n-1}} =
    [K_{n-1}^\times]$.  Now $\hat Y_{n-1} + \cdots + \hat Y_0
    \subset [K_{n-1}^\times]$ by construction, and $\hat
    Y_n^{H_{n-1}} = W = \hat Y_n \cap [K_{n-1}^\times] \subset
    [K_{n-1}^\times]$ from above.  Hence $\hat Y^{H_{n-1}} \subset
    [K_{n-1}^\times]$.  We also have the decomposition
    $[K_{n-1}^\times] = \epsilon(\tilde Y_{n-1}) + \hat Y_{n-2} +
    \cdots + \hat Y_0$.  Therefore it is sufficient to show that
    $\epsilon(\tilde Y_{n-1}) \subset \hat Y^{H_{n-1}}$.

    Because $\epsilon(\tilde Y_{n-1}) = \hat
    Y_{n-1} + P$ it is enough to show that $P\subset \hat
    Y_n^{H_{n-1}} + \hat Y_{n-2} + \cdots + \hat Y_0 =
    W + \hat Y_{n-2} + \cdots + \hat Y_0$.  But by the definition
    of projection, $P\subset W+\hat Y_{n-2} + \cdots + \hat Y_0$.
    Hence we conclude that $\hat Y^{H_{n-1}} = [K_{n-1}^\times]$,
    which is item (2) for $i=n-1$.

    For $i<n-1$, observe that since Theorems~\ref{th:nox} and
    \ref{th:x} hold in the case $n-1$, we have
    \begin{equation*}
        (\tilde Y_{n-1}+ \dots + \tilde Y_0)^ {H_i/H_{n-1}} =
        [K_i^\times]_{K_{n-1}}, \ i<n-1.
    \end{equation*}
    (If we are in the situation covered by Theorem~\ref{th:nox}
    then this statement is immediate. If we are in the situation
    covered by Theorem~\ref{th:x} we use the fact that $i(K_
    {n-1}/F)=-\infty$ and therefore the summand of $[K_i^
    \times]_{K_{n-1}}$ corresponding to the module generated
    by an exceptional element is trivial.)

    Again using Theorems~\ref{th:nox} and \ref{th:x} as well as the
    equality $\hat Y^{H_{n-1}}=[K_{n-1}^\times]$ and the fact that
    \begin{equation*}
        \epsilon\colon[K_{n-1}^\times]_{K_{n-1}}\to J
        \mbox{ with }\epsilon([K_{n-1}^\times]_{K_{n-1}
        })=[K_{n-1}^\times]
    \end{equation*}
    is an $\Fp[G]$-homomorphism, we obtain for each $i\in\{0,1,
    \dots,n-2\}$ that
    \begin{align*}
        \hat Y^{H_i} &=(\hat Y^{H_{n-1}})^{H_i/H_{n-1}}
        =[K_{n-1}^\times]^{H_i/H_{n-1}} \\
        &=(\epsilon(\tilde Y_{n-1}+\cdots+\tilde
        Y_0))^{H_i/H_{n-1}}\\ &=\epsilon([K_i^\times]_
        {K_{n-1}}) =[K_i^\times],
    \end{align*}
    as required.
\end{proof}

\begin{proof}[of Theorem~{\rm\ref{th:nox}}]
    The case $p=2$, $n=1$ was treated in Proposition~\ref{pr:noxn1}.
    For the remaining case of $\xi_p\notin F$ and $p>2$, we proceed
    by induction. The base case of $n=1$ is
    Proposition~\ref{pr:noxn1}. Assume then that $n>1$ and the
    Theorem holds for $n-1$.  By Proposition~\ref{pr:y} above, there
    exists an $\Fp[G]$-submodule $\hat Y= \oplus \hat Y_i \subset
    J$, where each $\hat Y_i$ is a direct sum of cyclic
    $\Fp[G]$-modules of dimension $p^i$, $[K_i^\times] = {\hat
    Y}^{H_i}$, $0\le i< n$, and $\hat Y_n^G = [N_{K/F}(K^\times)]$.
    Set $Y_i=\hat Y_i$ and $Y=\oplus \hat Y_i$. All that remains is
    to show that $J\subset Y$.

    We adapt the proof of the Inclusion Lemma~\eqref{le:incl} to
    show that $J\subset Y$, by induction on the socle series $J_i$
    of $J$. We first show that $J_{p^{n-1}}\subset Y$, as follows.
    Consider $Y$ and $J$ as $\Fp[H_{n-1}]$-modules. By the Fixed
    Submodule Lemma~\eqref{le:fixed}, $J^{H_{n-1}} =
    [K_{n-1}^\times]$, and we have already shown that
    $[K_{n-1}^\times] = Y^{H_{n-1}} \subset Y$, so
    $J^{H_{n-1}}\subset Y$.  Since $J^{H_{n-1}}$ is the kernel of
    $\sigma^{p^{n-1}}-1\equiv (\sigma-1)^{p^{n-1}}$, $J^{H_{n-1}}
    = J_{p^{n-1}}$.  Hence $J_{p^{n-1}} = [K_{n-1}^\times] \subset
    Y$.

    For the inductive step, assume that $J_i\subset Y$ for all
    $i<t$ for some $p^{n-1} < t \le p^n$, and let $[\gamma]\in
    J_{t}\setminus J_{t-1}$.  Hence $l(\gamma)=t$. Therefore
    $[\gamma]\not\in [K_{n-1}^\times]$, and by the Second Fixed
    Elements are Norms Lemma~\eqref{le:fixednew}, part (a), there
    exists $[\chi]\in J$ such that $[\gamma]^{(\sigma-1)^{t-1}} =
    [N_{K/F}(\chi)] \in Y_n^G$. Since $Y_n$ is a free
    $\Fp[G]$-module, there exists $[\chi']\in Y_n$ such that
    $[N_{K/F}(\chi')] = [\chi']^{(\sigma-1)^{p^n-1}} =
    [\chi]^{(\sigma-1)^{l(\chi)-1}}$. Set
    $[\gamma']=[\chi']^{(\sigma-1)^{p^n-t}} \in Y_n \subset Y$.
    Then $l(\gamma/\gamma')<t$.  By induction $[\gamma/\gamma']\in
    Y$, and since $[\gamma'] \in Y$, we obtain $[\gamma]\in Y$ as
    well.
\end{proof}

\section{Exceptional Elements}\label{se:exc}

Assume that $\xi_p\in F$ and, if $p=2$, then either $n>1$ or $-1
\in N_{K/F}(K^\times)$. Recall that in Proposition~\ref{pr:xn1} in
section~\ref{se:base} we proved that Theorem~\ref{th:x} holds for
extensions of degree $p$ and in Proposition~\ref{pr:xp2n2} we
proved that Theorem~\ref{th:x} holds in the case $p=2$ and
$n=2$. Assume then that Theorem~\ref{th:x} holds for extensions
of degree $p^s$ for $1\le s< n$.

In the next lemma we assume that $n\geq 2$ and, if $p=2$, that $n>
2$ as well. These conditions allow us to use Proposition~\ref{pr:y},
by which we assume that we have a submodule $\hat Y=\hat
Y_n\oplus\hat Y_{n-1}\oplus\cdots\oplus\hat Y_0$ of $J$ with
properties (1), (2) and (3) listed in Proposition~\ref{pr:y}.

\begin{lemma}\label{le:excdescent}
    Suppose $\delta\in K^\times$ satisfies $[N_{K/F}(\delta)]_F
    \neq [1]_F$ and $p^t+2\le l(\delta)\le p^{t+1}$, for some $t
    \in \{0, 1, \dots, n-2\}$. Then there exists $\delta'\in
    K^\times$ with $[N_{K/F}(\delta')]_F \neq [1]_F$ and
    $l(\delta') < l(\delta)$.
\end{lemma}

\begin{proof}
    Let $[\beta] = [\delta]^{ (\sigma-1)}$ and $[\gamma] =
    [\delta]^{(\sigma-1)^{l(\delta)-1}}$. Since $l(\beta) <
    p^{t+1}$, $[\beta]\in J^{H_{t+1}}$, and since $[\beta] \in
    J^{\sigma-1}$, $[N_{K/F}(\beta)]_F = [1]_F$.  By the Proper
    Subfield Lemma~\eqref{le:prop}, we have $[\beta]\in
    [K_{t+1}^\times]$.

    By Proposition~\ref{pr:y}, $[\beta] \in \hat Y^{H_{t+1}}$.
    Moreover, $p^{t} + 1\le l(\beta) < p^{t+1}$. Let
    \begin{align*}
        W &= \hat Y_{n}^{H_{t+1}} \oplus \hat Y_{n-1}^{H_{t+1}}
        \oplus \dots \oplus \hat Y_{t+1}^{H_{t+1}} \\ &= \hat
        Y_{n}^{H_{t+1}} \oplus \hat Y_{n-1}^{H_{t+1}} \oplus \dots
        \oplus \hat Y_{t+2}^{H_{t+1}} \oplus \hat Y_{t+1}.
    \end{align*}
    By Proposition~\ref{pr:y} and the Submodule-Subfield
    Lemma~\eqref{le:sub}, $W$ is a direct sum of cyclic
    $\Fp[G]$-modules of length $p^{t+1}$.

    Let $[\beta']$ be the component of $[\beta]$ in $W$.  Because
    $p^t+1\le l(\beta)$ and $(\sigma-1)^{p^t}$ is trivial on $\hat
    Y_{t} \oplus \dots \oplus \hat Y_0$, we see that $l(\beta) =
    l(\beta')$ and also
    \begin{equation*}
        [\gamma]=[\delta]^{(\sigma-1)^{l(\delta)-1}} =
        [\beta]^{(\sigma-1)^{l(\beta)-1}} = [\beta']^
    {(\sigma-1)^{l(\beta')-1}}.
    \end{equation*}
    Since $W$ is a direct sum of cyclic $\Fp[G]$-modules of length
    $p^{t+1}$ and contains $[\beta']$, of length strictly less
    than $p^{t+1}$, $[\beta']$ lies in the image of $(\sigma-1)$
    on $W$. Hence there exists $[\alpha']\in W$ such that
    $[\alpha']^{(\sigma-1)} = [\beta']$. Therefore $l(\alpha')
    =l(\delta)$ and
    \begin{equation*}
        [\gamma]=[\delta]^{(\sigma-1)^{l(\delta)
    -1}} = [\alpha']^{(\sigma-1)^{l(\delta)-1}}.
    \end{equation*}
    Moreover, by Proposition~\ref{pr:y}, $W\subset \hat Y^{H_{t+1}}
    = [K_{t+1}^\times] \subset [K_{n-1}^\times]$ and therefore we
    have $[N_{K/F}(\alpha')]_F = [1]_F$. Now set $\delta' =
    \delta/\alpha'$.  Then $[\delta']^{(\sigma-1)^{l(\delta)-1}}
    =[1]$ so that $l(\delta') < l(\delta)$, and
    $[N_{K/F}(\delta')]_F = [N_{K/F}(\delta)]_F \neq [1]$.
\end{proof}

\begin{proposition}\label{pr:ldelta}
    Suppose that $\xi_p\in F$ and, if $p=2$, that $n>1$ or $-1\in
    N_{K/F}(K^\times)$. Then $m<n$ and, for any exceptional element
    $\delta$, $l(\delta) = p^{m}+1$. Moreover, this length is the
    minimal $l(z)$ for all $z\in K^\times$ with $[N_{K/F}(z)]_F \neq
    [1]_F$.
\end{proposition}

Observe that the proposition implies that for any exceptional
element $\delta$, $l(\delta)<p^n$.  (Indeed $p^m+1\le p^{n-1}+1$ and
$p^{n-1}+1\le p^n$ unless $p=2$ and $n=1$.  If $p=2$, $n=1$, and
$-1\in N_{K/F}(K^\times)$, then let $-1=N_{K/F}(\theta)$, where
$\theta\in K^\times$. Observe that $\delta = \sqrt{a} \theta$
satisfies $[N_{K/F}(\delta)]_F = [a]_F$ and $[\delta]^ {(\sigma-1)}
= [1]$. Hence $l(\delta)<2$.) By the Norm Lemma~\eqref{le:norm},
then $[N_{K/F}(\delta)]_F = [a]_F^s$, and by definition of
exceptional element, $s\not\equiv 0\bmod p$. By choosing an
appropriate power of $\delta$, we have there exists an exceptional
element $\delta$ with $[N_{K/F}(\delta)] _F = [a]_F$.

\begin{proof}
    We first prove that $m<n$.  Assume first that $p>2$ or $p=2$
    and $n>1$.  Consider $\delta = \root{p}\of{a_{n-1}}$.  We
    observed in the proof of Proposition \ref{pr:excexist} that
    $N_{K/F}(\delta) = a_0 = a$.  Now $\delta^\sigma =
    \root{p}\of{a_{n-1}^\sigma}$ for a suitable $p$th root of
    unity.  Because $K/K_{n-1}$ is Galois we see from Kummer
    theory that $a_{n-1}^\sigma = a_{n-1} k_{n-1}^p$ for some
    $k_{n-1} \in K^\times_{n-1}$.  Hence $\delta^{\sigma-1} \in
    K_{n-1}^\times$, and therefore $m \leq n-1$, as required.

    Proposition~\ref{pr:ldelta} was established for the case $n=1$,
    and if $p=2$, also for $n=2$, in Propositions~\ref{pr:xn1} and
    \ref{pr:xp2n2}.  Therefore we now assume that $n\ge 2$ and if
    $p=2$ then $n>2$ as well.

%
    Now let $\delta$ be an arbitrary exceptional element.
    Clearly $[\delta]\neq [1]$ since $[N_{K/F}(\delta)]_F \neq
    [1]_F$; hence $l(\delta)\ge 1$. If $m=-\infty$, then
    $[\delta]^{\sigma-1} = [1]$ so that $l(\delta)\le 1$ and
    because of our convention $p^{-\infty}=0$ we are done.

    Hence assume that $m\geq 0$. Then set $[\beta] :=
    [\delta]^{\sigma-1} \in [K_m^\times]$. Also
    \begin{equation*}
        [N_{K_m/F}(\beta)] = [\beta]^{(\sigma-1)^{p^m-1}} \in
        [F^\times].
    \end{equation*}
    Therefore $l(\delta)\leq 1+ (p^m-1) + 1 = p^m+1$.

    Now suppose that $[z]\in J$ satisfies $[N_{K/F}(z)]_F \neq
    [1]_F$ and $l(z)$ is minimal among all such $z$.  Since
    $[N_{K/F}(\delta)]_F \neq [1]_F$ and $\delta$ above has
    $l(\delta) \le p^{m}+1$, we see that $l(z) \le p^{m}+1$. Now
    suppose, contrary to our statement, that $l(z) < p^{m} + 1$. If
    $m=0$ then $l(z) = 1$ and hence $[z]^{(\sigma-1)} \in
    [K_{-\infty}^\times]$, contradicting the minimality of $m$.
    Otherwise $m\ge 1$ and repeated application of
    Lemma~\ref{le:excdescent} yields $\delta'\in K^\times$ such that
    $[N_{K/F}(\delta')]_F \neq [1]_F$ and $l(\delta') \le p^{m-1} +
    1$. (Observe that we can indeed apply Lemma~\ref{le:excdescent},
    since $l(\delta') \leq l(z) \leq p^{m} \leq p^{(n-2)+1}$, where
    the last inequality holds since $m <n$.)

    Let $[\beta']=[\delta']^{(\sigma-1)}$.  Then $l(\beta') \le
    p^{m-1}$ so that $[\beta']\in J^{H_{m-1}}$, and since
    $[\beta']$ is in the image of $(\sigma-1)$,
    $[N_{K/F}(\beta')]_F = [1]_F$. By the Proper Subfield
    Lemma~\eqref{le:prop}, we see that $[\beta'] \in
    [K_{m-1}^\times]$.  Hence $[N_{K/F}(\delta')]_F \neq [1]_F$,
    and $[\delta']^{\sigma-1} \in [K_{m-1}^\times]$, contradicting
    the minimality of $m$.  Therefore $l(\delta)=p^{m}+1$.
\end{proof}

Now assume that $\xi_p\in F$ and, if $p=2$, then $n\ge 2$.

\begin{proposition}\label{pr:excexc}
    If $\delta$ is an exceptional element of $K/F$, then
    $\delta$ is an exceptional element of $K/K_i$ for $0\le i
    < n$ if $p>2$ and for $0\le i < n-1$ if $p=2$.
\end{proposition}

\begin{proof}
    Since $K_0=F$, the proposition is clear for $i=0$.  We
    therefore assume that $i>0$.

    If $\delta$ is an exceptional element of $K/F$, then
    Proposition~\ref{pr:ldelta} tells us that $l(\delta) = p^m+1$
    for $m<n$. If $p>2$, then for each $i\in\{0,1,\dots,n-1\}$ we
    have
    \begin{equation*}
        l(\delta) = p^m+1 \leq p^{n-1}+1 \leq p^n-p^{n-1} \leq
        p^n-p^i.
    \end{equation*}
    If $p=2$ and $n\geq 2$, then similarly for each
    $i\in\{0,1,\dots,n-2\}$ we have
    \begin{equation*}
        l(\delta) = 2^m+1 \leq 2^{n-1}+1 \leq 2^n-2^{n-2} \leq
        2^n-2^i.
    \end{equation*}
    Since the Norm Lemma~\eqref{le:norm} gives
    $[N_{K/F}(\delta)]_F \neq [1]_F$ if and only if
    $[N_{K/K_i}(\delta)]_{K_i} \neq [1]_{K_i}$, it follows
    $[N_{K/K_i}(\delta)]_{K_i}\neq [1]_{K_i}$.

    Let $\tau = \sigma^{p^i}$.  Then $(\tau - 1) \equiv
    (\sigma-1)^{p^i}$ on $J$, and so $[\delta]^{(\sigma-1)}\in
    [K_m^\times]$ implies that $[\delta]^{(\tau-1)} \in
    [K_m^\times]$. Now we define intermediate fields
    $\{K'_{-\infty},K'_0, \dots, K'_{n-i}\}$ of $K/K_i$ by $K'_j
    := K_{j+i}$.

    First consider the case $m<i$.  We have $l(\delta) = p^m+1$, and
    so then $l(\delta) \le p^i$.  Hence $[\delta]^{\tau-1} =
    [\delta]^{(\sigma-1)^{p^i}} = [1]$. Since we have shown that
    $\delta$ satisfies $[N_{K/K_i}(\delta)]_{K_i} \neq [1]_{K_i}$,
    $\delta$ is an exceptional element of $K/K_i$ with
    $i(K/K_i)=-\infty$.

    Now consider the case $m\ge i$.  In this case we have shown
    that $\delta$ satisfies $[N_{K/K_i}(\delta)]_{K_i} \neq
    [1]_{K_i}$ and $[\delta]^{\tau-1} \in [K_{m-i}^{'\times}]$.
    All that remains is to show that no $\delta'\in K^\times$
    exists with $[N_{K/K_i} (\delta')]_{K_i}\neq [1]_{K_i}$ and
    $[\delta']^{\tau-1} \in [K_{j}^{'\times}]$ for $j < m-i$.
    Suppose such a $\delta'$ exists. We may assume that this
    $\delta'$ has a minimal length among all elements $z$ with
    $[N_{K/K_i}(z)]_{K_i} \neq [1]_{K_i}$. By the remark made
    after Proposition~\ref{pr:ldelta} we see that we may further
    assume that $[N_{K/K_i}(\delta')]_{K_i} = [a_i]_{K_i}$.
    Therefore $[N_{K/F}(\delta')]_F = [N_{K_i/F}(a_i)]_F =
    [a]_F\neq [1]_F$.

    If $j=-\infty$, then since $(\tau-1)\equiv (\sigma-1)^{p^i}$,
    we obtain $l(\delta') \le p^i \le p^m$.  On the other hand, if
    $j\ge 0$ then $m>i$. Moreover, since $(\tau-1)\equiv
    (\sigma-1)^{p^i}$ and $[N_{K'_j/F}(\gamma)] =
    [\gamma]^{(\sigma-1)^{p^{i+j}-1}}$ for $[\gamma]\in
    [K_j^{'\times}]$, we have
    \begin{equation*}
        l(\delta')\le p^i + (p^{i+j}-1) + 1 \le p^m.
    \end{equation*}
    In either case, this violates the condition of
    Proposition~\ref{pr:ldelta}, since then $l(\delta)$ is not
    minimal among lengths $l(\delta')$ for $[N_{K/F}(\delta')]_F
    \neq [1]_F$.
\end{proof}

\section{Proof of Theorem~2}\label{se:x}

We first adapt the proof of Theorem~\ref{th:nox} to prove the
following analogue for the case of Theorem~\ref{th:x}.  We assume
here that Theorem~\ref{th:x} holds for $n-1$ and, if $p=2$, then
$n>2$.

\begin{proposition}\label{pr:xhaty}
    Let $\xi_p\in F$, $n\ge 2$, and $\delta\in K^\times$ be any
    exceptional element of $K/F$. Then the $\Fp[G]$-module $J$
    decomposes as
    \begin{equation*}
        J=X + \hat Y, \quad \hat Y = \hat Y_n \oplus \hat Y_{n-1}
        \oplus \dots \oplus \hat Y_0,
    \end{equation*}
    where
    \begin{enumerate}
        \item $X$ is the cyclic $\Fp[G]$-module generated by
        $[\delta]$;
        \item for each $i=1, 2, \dots, n$, $\hat Y_i$ is a direct
        sum of cyclic $\Fp[G]$-modules of dimension $p^i$;
        \item $[K_i^\times] = \hat Y^{H_i}$ for $0\le i < n$; and
        \item $\hat Y_n^G = [N_{K/F}(K^\times)]$.
    \end{enumerate}
\end{proposition}

\begin{proof}
    By Proposition~\ref{pr:y}, there exists an $\Fp[G]$-submodule
    $\hat Y= \oplus \hat Y_i \subset J$, where each $\hat Y_i$ is
    a direct sum of cyclic $\Fp[G]$-modules of dimension $p^i$,
    $[K_i^\times] = {\hat Y}^{H_i}$, $0\le i< n$, and $\hat Y_n^G
    = [N_{K/F}(K^\times)]$.  Let $X$ be defined as in the
    statement of the Theorem and set $\hat J = X + \hat Y$. We
    have $\hat J$ is an $\Fp[G]$-submodule of $J$.

    Assume first that $p>2$.  Consider $\hat J$ and $J$ as
    $\Fp[H_{n-1}]$-modules.  By Proposition~\ref{pr:excexc},
    $\delta$ is exceptional for $K/K_{n-1}$ and so by
    Proposition~\ref{pr:xn1} (which is just Theorem~\ref{th:x}
    in case $n=1$), $J$ decomposes as $\bar X \oplus
    \bar Y_1 \oplus \bar Y_0$, where $\bar X\subset X$ is the
    $\Fp[H_{n-1}]$-submodule generated by $[\delta]$, $\bar
    Y_0\subset J^{H_{n-1}}$, and $\bar Y_1^{H_{n-1}} + \bar Y_0
    \subset [K_{n-1}^\times]$ (by Theorem~\ref{th:x}, case
    $n=1$, part (3)). Hence $J^{H_{n-1}} \subset X +
    [K_{n-1}^\times] \subset \hat J$. (Here we use the fact
    that $\bar X\subset X$ and $\bar Y_0^{H_{n-1}}=\bar Y_0$.)

    Now suppose that $[\Gamma]\in J\setminus
    (X+[K_{n-1}^\times])$. Our goal is to show that $[\Gamma] =
    [\theta] + [\gamma]$ with $[\theta]\in X$ and
    $[\gamma]^{(\sigma-1)^{l(\gamma)-1}}\in \hat Y_n^G$. Then,
    with this result in hand, we will adapt the proof of the
    Inclusion Lemma~\eqref{le:incl} to show that $J\subset \hat
    J$.

    Write $l_H(\Gamma)$ for the length of the cyclic
    $\Fp[H_{n-1}]$-submodule of $J$ generated by $[\Gamma]$.  Since
    $[\Gamma]\not\in J^{H_{n-1}}$, we find $l_H(\Gamma) \ge 2$.

    If $l_H(\Gamma)=2$ and $\Gamma$ is exceptional, we find $\gamma$
    and $\theta$ as follows. By Proposition~\ref{pr:excexc},
    $\delta$ and $\Gamma$ are exceptional elements for $K/K_{n-1}$.
    Since $[N_{K/K_{n-1}}(\delta)]_{K_{n-1}} \neq [1]_{K_{n-1}}$, we
    see that for a suitable power $s\in \Z$, $[N_{K/K_{n-1}}
    (\Gamma)]_{K_{n-1}} = [N_{K/K_{n-1}} (\delta)]_{K_{n-1}}^s$. Set
    $\theta=\delta^s$ and $\gamma = \Gamma/\theta$. Then
    $[N_{K/K_{n-1}}(\gamma)]_{K_{n-1}} = [1]_{K_{n-1}}$ and so
    $[N_{K/F}(\gamma)]_F = [1]_F$.  Moreover, $l(\gamma)> p^{n-1}$
    since otherwise $[\gamma]\in J^{H_{n-1}}$ and by the Exact
    Sequence Lemma~\eqref{le:exact}(ii) we would have $[\gamma] \in
    [K_{n-1}^\times]$, contradicting our assumption on $\Gamma$.
    Thus we have $l(\gamma),l(\Gamma)> p^{n-1}$. Also since the
    maximum length of the elements in $X$ is at most $p^{n-1}+1$ by
    Proposition~\ref{pr:ldelta}, we have $l(\theta)\leq p^{n-1}+1$.
    Now if $l(\Gamma) > l(\theta)$ then $l(\gamma) =
    l(\Gamma/\theta)=l(\Gamma)$, and if $l(\Gamma) = l(\theta)$ then
    $l(\Gamma)=l(\theta)=p^{n-1} +1$ and $p^{n-1} < l(\gamma)\leq
    p^{n-1}+1$, showing that in this case as well $l(\Gamma) =
    l(\gamma)$. Thus in all cases $l(\Gamma) = l(\gamma)$ and
    therefore also $l_H(\gamma)=l_H( \Gamma)$.

    Otherwise, let $\theta=1$ and $\gamma=\Gamma$.  Clearly
    $l(\gamma)=l(\Gamma)$ and $l_H(\gamma) = l_H(\Gamma)$.

    In either case, our choice of $\gamma$ is made in order to make
    sure that we have either $l_H(\gamma)\geq 3$ or both
    $l_H(\gamma)=2$ and $[N_{K/F}(\gamma)]_F = [1]_F$. These are the
    necessary hypotheses to apply the Second Fixed Elements are
    Norms Lemma~\eqref{le:fixednew}, part (a), by which we obtain
    that there exists $[\chi]\in J$ such that
    $[\gamma]^{(\sigma-1)^{l(\gamma)-1}} = [N_{K/F}(\chi)] \in \hat
    Y_n$.  Hence we have shown that for all $[\Gamma]\in J\setminus
    (X+[K_{n-1}^\times])$, we have $[\Gamma]=[\theta]+[\gamma]$ with
    $[\theta]\in X$, $[\gamma]^{(\sigma-1)^{l(\gamma)-1}}\in \hat
    Y_n^G$.

    Now we adapt the proof of the Inclusion Lemma~\eqref{le:incl}
    to show that $J\subset \hat J$, by induction on the socle
    series $J_i$ of $J$. Since $\sigma^{p^{n-1}}-1\equiv
    (\sigma-1)^{p^{n-1}}$, $J^{H_{n-1}}=J_{p^{n-1}}$.  Hence
    $J_{p^{n-1}}\subset \hat J$ and our base case for the
    induction is $J_{p^{n-1}}$.

    For the inductive step, assume that $J_i\subset \hat J$ for all
    $i<t$ for some $p^{n-1} < t \le p^n$, and let $[\Gamma]\in
    J_{t}\setminus J_{t-1}$.  Then $l(\Gamma)=t$.  If $[\Gamma] \in
    X + [K_{n-1}^\times]$, we have already shown that $[\Gamma] \in
    \hat J$.  Therefore we assume that this is not the case.

    By our result above, we may write $[\Gamma] = [\theta]+[\gamma]$
    with $[\theta]\in X$ and
    \begin{equation*}
        [\gamma]^{(\sigma-1)^{l(\gamma)-1}} = [N_{K/F}(\chi)]
        \in \hat Y_n^G
    \end{equation*}
    for some $[\chi]\in J$. Moreover, as we have shown, we may
    assume that $l(\gamma) = l(\Gamma)$. To show that $[\Gamma]\in
    \hat J$, it is enough to show that $[\gamma]\in \hat J$. Since
    $\hat Y_n$ is a free $\Fp[G]$-module, there exists $[\chi']\in
    \hat Y_n$ such that $[N_{K/F}(\chi')] =
    [\chi']^{(\sigma-1)^{p^n-1}} =
    [\chi]^{(\sigma-1)^{l(\chi)-1}}$. Set $[\gamma'] =
    [\chi']^{(\sigma-1)^{p^n-t}} \in \hat Y_n \subset \hat J$.
    Then $l(\gamma/\gamma')<t$.  By induction $[\gamma/\gamma']\in
    \hat J$, and since $[\gamma'] \in \hat J$, $[\gamma]\in \hat
    J$ as well. Hence our induction is complete.

    The case $p=2$ follows similarly with the following
    modifications. Replace $H:=H_{n-1}$ with $H:=H_{n-2}$ and
    $K_{n-1}$ with $K_{n-2}$. Thus we consider $\hat J$ and $J$ as
    $\Ft[H]$-modules. By Proposition~\ref{pr:xp2n2} (our theorem
    in the base case $p=2$ and $n=2$) and by
    Proposition~\ref{pr:excexc} we may write
    \begin{equation*}
        J = \bar X \oplus \bar Y_0 \oplus \bar Y_1 \oplus \bar
        Y_2,
    \end{equation*}
    where $\bar X\subset X$ is the cyclic $\Ft[H]$-module
    generated by $[\delta]$ and for $i=0, 1, 2$, the summand $\bar
    Y_i$ is a direct sum of cyclic $\Ft[H]$-modules of dimension
    $2^i$. By Proposition~\ref{pr:xp2n2} we also have
    \begin{align*}
        J^{H_{n-1}} &\subset \bar X^{H_{n-1}} \oplus (\bar Y_0
        \oplus \bar Y_1 \oplus \bar Y_2)^{H_{n-1}} \\ &\subset
        X\oplus [K_{n-1}^\times]\subset\hat J.
    \end{align*}

    Now suppose that $[\Gamma]\in J\setminus (X+[K_{n-1}^\times])$.
    Again we want to show that $[\Gamma] = [\theta]+[\gamma]$ with
    $[\theta]\in X$ and $[\gamma]^{(\sigma-1)^{l(\gamma)-1}} \in
    \hat Y_n^G$. We have $[\Gamma]\notin J^{H_{n-1}}$ and so
    $l_H(\Gamma) \geq 3$.

    If $l_H(\Gamma)=3$ then
    \begin{equation*}
        [N_{K/K_{n-2}}(\Gamma)]_{K_{n-2}} = [a_{n-2}]_{K_{n-2}}^s
    \end{equation*}
    for some $s \in \Z$. Set $\theta=\delta^s$ and
    $\gamma=\Gamma/\theta$. Then
    \begin{equation*}
        [N_{K/K_{n-2}}(\gamma)]_{K_{n-2}} = [1]_{K_{n-2}},
    \end{equation*}
    whence $[N_{K/F}(\gamma)]_F = [1]_F$.
    Also $[\gamma]\notin X + [K_{n-1}^\times]$ and therefore
    $l(\gamma)>2^{n-1}$. On the other hand, $l(\theta)\leq
    2^{n-1}+1$ by Proposition~\ref{pr:ldelta}. Hence we see again
    that $l(\gamma)=l(\Gamma)$ and in particular $l_H(\gamma)\geq
    3$.

    Otherwise, if $l_H(\Gamma)=4$ then let $\theta=1$ and
    $\gamma=\Gamma$. Clearly $l(\gamma)=l(\Gamma)$ and
    $l_H(\gamma) = l_H(\Gamma)$.

    In either case, we see from the Second Fixed Elements are Norms
    Lemma~\eqref{le:fixednew}, part (b), that there exists
    $[\chi] \in J$ such that $[N_{K/F}(\chi)] =
    [\gamma]^{(\sigma-1)^{l(\gamma)-1}}$.

    From now on the proof that $\hat J=J$ in the case $p=2$ is
    identical with the proof above for the case $p>2$.
\end{proof}

\begin{proof}[of Theorem~{\rm\ref{th:x}}]
    The case $n=1$ is Proposition~\ref{pr:xn1}. The case $p=2$ and
    $n=2$ was established in Proposition~\ref{pr:xp2n2}. We proceed
    by induction on $n$. Assume therefore that the Theorem holds for
    $n-1$, and assume that $n>1$ if $p>2$ and $n>2$ if $p=2$.  By
    Proposition~\ref{pr:xhaty}, we write $J = X+\hat Y$, $\hat
    Y=\oplus \hat Y_i$, where $\hat Y_i$ is a direct sum of cyclic
    $\Fp[G]$-modules of dimension $p^i$ and for $i<n$, $[K_i^\times]
    = \hat Y^{H_i}$.

    We define the $Y_i$ and $Y$ as follows.  When $m = -\infty$,
    set $Y_i = \hat Y_i$ and $Y = \sum Y_i$.  Now suppose that $m
    \geq 0$, and let $[\beta]\in [K_m^\times]$ satisfy
    $[\beta]=[\delta]^{(\sigma-1)}$. By Proposition~\ref{pr:ldelta},
    we see that $l(\beta)=p^m$ and so the cyclic $\Fp[G]$-submodule
    $M_\beta$ generated by $[\beta]$ is a free $\Fp[G/H_m]$-submodule.
    Moreover, we have already established that $[K_m^\times]\subset
    \hat Y$, so $M_\beta \subset \hat Y$.

    Now let $[\gamma]=[\beta]^{(\sigma-1)^{p^m-1}} \in \hat Y^G$.
    Suppose that $[\gamma]\in W := \hat Y_{m+1}+\dots+\hat Y_n$.
    Assume first that if $p=2$ then $m>0$. Then since $W^G$ is in
    the image of $(\sigma-1)^{p^m+1}$ on $W$, there exists
    $[\alpha]\in W$ such that $[\alpha]^{(\sigma-1)^{p^m+1}} =
    [\gamma]$.  Hence $[\beta'] = [\alpha]^{(\sigma-1)}$, being in
    the image of $(\sigma-1)$, satisfies $[N_{K/F}(\beta')]_F =
    [1]_F$, while $l(\beta') = p^m+1 = l(\delta)$ and
    $[\beta']^{(\sigma-1)^{p^m}} = [\delta]^{(\sigma-1)^{p^m}}$.
    Hence $[N_{K/F}(\delta/\beta')]_F \neq [1]_F$ and
    $l(\delta/\beta') \le p^m$.  But this contradicts the minimality
    of $l(\delta)$ among lengths $l(z)$ with $[N_{K/F}(z)]_F \neq
    [1]_F$, a contradiction. Hence $[\gamma] \not\in W$.

    Now if $p=2$ and $m=0$ we have $[\beta] = [\gamma]$.  Suppose
    $[\beta] \in W$. Then there exists $[\nu] \in W := \hat Y_1
    +\dots+ \hat Y_n$ such that $[\nu]^{\sigma-1} = [\beta]$. Hence
    $[\nu] = [\delta\nu] [\delta]$ and $[\delta\nu] \in J^G$.  Since
    $l(\delta\nu) < l(\delta)$, Proposition~\ref{pr:ldelta} tells us
    that $[N_{K/F}(\delta\nu)]_F = [1]_F$. Hence $\nu$ is
    exceptional. On the other hand, $[\nu] \in \hat Y^{H_1} \subset
    [K_1^\times]$.  But since $n>1$, $N_{K/F}(K_1^\times)\subset
    F^{\times p}$, so $[N_{K/F}(\nu)]_F = [1]_F$, contradicting the
    exceptionality of $\nu$.  Thus again $[\gamma] = [\beta] \not\in
    W$.

    If $m=0$ then $[\gamma]=[\beta]$ lies in $\hat Y_0\oplus\hat
    Y_1\oplus \cdots \oplus \hat Y_n$.  If $m\ge 1$ then $[\gamma]$
    is similarly in $\hat Y_m \oplus \cdots \oplus \hat Y_n$ since
    $[\gamma]$ is in the image of $(\sigma-1)^{p^m-1}$ on $\hat Y$.
    Let $[\gamma']$ be the component of $[\gamma]$ in $\hat Y_m$. By
    the previous two paragraphs, $[\gamma']\neq [1]$.  Now since
    $[\gamma']$ lies in $\hat Y_m^G$, we have $[\gamma']$ is in the
    image of $(\sigma-1)^{p^m-1}$ on $\hat Y_m$. Now let
    $[\beta]_{(m)}\in \hat Y_m$ be the projection of $[\beta]$ into
    $\hat Y_m$. (Since $M_\beta\subset\hat Y$ this projection is
    well defined.) Moreover since $[\gamma'] =
    [\beta]_{(m)}^{(\sigma-1)^{p^m-1}} \neq [1]$ we see that
    $[\beta]_{(m)}$ generates a cyclic $\Fp[G]$-submodule
    $M_{[\beta]_{(m)}}$ of $\hat Y_m$ which is a free
    $\Fp[G/H_m]$-submodule of $\hat Y_m$. By the Free Complement
    Lemma~\eqref{le:free}, there exists a free
    $\Fp[G/H_m]$-complement $Y_m$ of $M_{[\beta]_{(m)}}$ in $\hat
    Y_m$. Having defined $Y_m$, we set all other $Y_i=\hat Y_i$,
    $i\neq m$, and $Y=\sum Y_i$.

    Since the $\hat Y_i$ are all independent, the $Y_i$ are
    independent. Assume now that $m\ge 0$. Then $X + \sum Y_i = X +
    \sum \hat Y_i$, because clearly $X + \sum Y_i \subset X + \sum
    \hat Y_i$, and $\hat Y_m \subset X + \sum Y_i$ follows from our
    construction of $Y_m$. Hence we have $X+Y =X+\hat Y=J$. Because
    in the case $m=-\infty$ we set $Y_i= \hat Y_i$ for all
    $i\in\{0,1,\dots,n\}$ we see that $J=X+\hat Y= X+Y$ as well. To
    show that the sum is direct, consider first the case $m =
    -\infty$. Here $X^G = X$, and by the Fixed Submodule
    Lemma~\eqref{le:fixed} and the equality $Y^G=[F^\times]$ we have
    $X\cap Y^G = \{0\}$.  Hence $X$ and $Y$ are independent. When $m
    \geq 0$, $X^G$ is generated by $[\gamma]$, which by construction
    of $Y$ satisfies $[\gamma]\not\in Y^G$. Using the Exclusion
    Lemma~\eqref{le:excl}, we have that $X$ and $Y$ are independent.

    We now show that $X^{(\sigma-1)}\oplus Y^{H_i}=\hat Y^{H_i}$ for
    $i\geq m$. (Here $X^{(\sigma-1)}$ means the image of $X$ under
    $(\sigma-1)$.) First observe that since $X$ and $Y$ are
    independent we indeed have $X^{(\sigma-1)}+Y^{H_i}=X^{(\sigma
    -1)}\oplus Y^{H_i}$. If $m=-\infty$ then $X^{(\sigma-1)}=\{0\}$
    and the equality $X^{(\sigma-1)}\oplus Y^{H_i}=\hat Y^{H_i}$ is
    a trivial statement. Assume now that $m\geq 0$. We have
    \begin{equation*}
        X^{(\sigma-1)}\subset[K_m^\times]\subset [K_i^
        \times]\subset\hat Y^{H_i}.
    \end{equation*}
    Hence for $m\le i\le n-1$
    \begin{equation*}
        X^{(\sigma-1)}\oplus Y^{H_i}\subset\hat Y^{H_i}.
    \end{equation*}
    To obtain the reverse inclusion, observe that $Y_i=\hat Y_i$,
    $i\neq m$, and $\hat Y_m\subset X^{(\sigma-1)}+Y_m$ by our
    construction of $Y_m$. Finally since $\hat Y_m^{H_i}=\hat Y
    _m$ and $Y_m^{H_i}=Y_m$ we see that also $\hat Y^{H_i}\subset
    X^{(\sigma-1)}\oplus Y^{H_i}$. Thus we indeed have the desired
    equality
    \begin{equation*}
        X^{(\sigma-1)}\oplus Y^{H_i}=\hat Y^{H_i}=[K_i^
        \times],
    \end{equation*}
    for each $i\in\{m,\dots,n-1\}$ if $m\geq 0$ and for each
    $i\in\{0, 1, \dots, n-1\}$ if $m=-\infty$. For $i<m$, observe
    that $X$ is cyclic of length $p^m+1$ and $(\sigma^{p^i}-1)\equiv
    (\sigma-1)^{p^i}$ on $J$. Then from the
    Submodule-Subfield Lemma~\eqref{le:sub} applied to the field
    extension $K_m/K_i$, we obtain
    \begin{equation*}
        X^{H_i} = X^{(\sigma-1)(\sigma^{p^i}-1)^{p^{m-i}-1}}.
    \end{equation*}
    Then, since $[K_i^\times]=\hat Y^{H_i} = (X^{(\sigma-1)})^{H_i}
    \oplus Y^{H_i}$ for all $i\le m$, we are done.
\end{proof}

\section{Proofs of Corollaries 1 and 2}

\begin{proof}[of Corollary~{\rm\ref{co:nox}}]
    Recall that if $M$ is a cyclic $\Fp[G]$-module of dimension
    $l$, then the $l+1$ submodules of $M$ are cyclic, given by
    $(\sigma-1)^{i}M$, $i=0, 1, \dots, l$, and have annihilators
    $\langle (\sigma-1)^{l-i}\rangle\subset \Fp[G]$, respectively.
    By Theorem~\ref{th:nox}, $[K_i^\times] = J^{H_i} = \oplus
    Y_i^{H_i}$, and $Y_j$ is a direct sum of cyclic
    $\Fp[G]$-modules of dimension $p^j$.

    Now $H_i=\langle \sigma^{p^i}\rangle$ and $(\sigma^{p_i}-1)
    \equiv (\sigma-1)^{p^i}$ on $J$.  When $j<i$, observe that $Y_j$
    is a direct sum of cyclic $\Fp[G]$-modules of dimension
    $p^j<p^i$, and so $Y_j=Y_j^{H_i}$.  When $j\ge i$, the submodule
    $Y_j^{H_i}$ is given by $Y_j^{(\sigma-1)^{p^j-p^i}}$.

    On $[K_i^\times]$, $N_{K_i/F} \equiv (\sigma-1)^{p^i-1}$. For
    $j<i$, since $Y_j^{H_i}$ is a direct sum of cyclic
    $\Fp[G]$-modules of dimension $p^j<p^i$, $N_{K_i/F}$
    annihilates $Y_j^{H_i}$. For $j\ge i$, $Y_j^{H_i}$ is a direct
    sum of cyclic $\Fp[G]$-modules of dimension $p^i$ and so
    applying $N_{K_i/F}$ to $Y_j^{H_i}$ yields
    $Y_j^{(\sigma-1)^{p^j-1}}=Y_j^G$. Hence we have the first
    statement.

    Now a cyclic $\Fp[G]$-module of dimension $p^i$ is a free
    $\Fp[G/H_i]$-module on one generator, and for direct sums $M$
    of such modules
    \begin{equation*}
        \rank_{\Fp[G/H_i]} M = \dim_{\Fp} M^G.
    \end{equation*}
    Observe that $Y_j$, $j<n$, is a direct sum of cyclic
    $\Fp[G]$-modules of dimension $p^j<p^n$.  Applying $N_{K/F}$
    to $J$, then, we see that $Y_n^G = [N_{K/F}(K^\times)]$.
    Moreover, with a descending induction we see that $Y_i^G$ is a
    complement of $[N_{K_{i+1}/F} (K_{i+1}^\times)]$ in
    $[N_{K_i/F}(K_i^\times)]$. Hence we have the second statement.
\end{proof}

\begin{proof}[of Corollary~{\rm\ref{co:x}}]
    We begin as in the previous proof.  If $m=-\infty$, in fact,
    then the previous proof carries over without modification.
    Hence we assume that $m\ge 0$.

    By Theorem~\ref{th:x}, $[K_i^\times] = X' \oplus Y_i^{H_i}$,
    where $X'$ is a cyclic $\Fp[G]$-module of dimension $p^i$ if
    $i\le m$ and of dimension $p^m$ if $i\ge m$. As in the
    previous proof, $Y_j^{H_i}=Y_j$ for $j< i$ and $Y_j^{H_i}$ for
    $j\ge i$ is a direct sum of cyclic $\Fp[G]$-modules of
    dimension $p^i$. Similarly, $N_{K_i/F}$ annihilates
    $Y_j^{H_i}$, $j<i$, and yields $Y_j^G$ when $j\ge i$. Applying
    $N_{K_i/F}$ annihilates $X'$ when $m<i$ and otherwise yields
    $(X')^G=X^G$.  Hence we have the statements locating
    $[N_{K_i/F}(K_i^\times)]$.

    For the statements establishing ranks, we proceed as in the
    previous proof. Observe that since $X\cap [K_m^\times] =
    X^{\sigma-1}$ is a cyclic $\Fp[G]$-submodule of dimension $p^m$,
    we obtain $X^G\subset [N_{K_m/F}(K_m^\times)]$. If $m\neq n-1$,
    then since $X\cap [K_{m+1}^\times] = X^{\sigma-1}$ is a cyclic
    $\Fp[G]$-submodule of dimension $p^m<p^{m+1}$, we see that
    $X^G\cap [N_{K_{m+1}/F}(K_{m+1}^\times)]=\{0\}$.  If $m=n-1$
    then $X\cap [K_{m+1}^\times]=X$ is a cyclic $\Fp[G]$-submodule
    of dimension $p^m+1$, which is annihilated by $N_{K_{m+1}/F}$
    unless $p^m+1=p^{m+1}=p^n$---that is, $p=2$, $m=0$, $n=1$. But
    this latter case violates the hypothesis of Theorem~\ref{th:x}.
    Hence $X^G\subset [N_{K_m/F}(K_m^\times)] \setminus
    [N_{K_{m+1}/F}(K_{m+1}^\times)]$ under our hypotheses.

    Again, since $Y_j$, $j<n$, and $X$ are direct sums of cyclic
    $\Fp[G]$-submodules of dimension less than $p^n$, applying
    $N_{K/F}$ to $J$ yields $Y_n^G=[N_{K/F}(K^\times)]$. A
    descending induction yields that $Y_i^G$, $m<i< n$, is a
    complement of $[N_{K_{i+1}/F}(K_{i+1}^\times)]$ in
    $[N_{K_i/F}(K_i^\times)]$. But $Y_m^G$ is a complement of
    $[N_{K_{m+1}/F}(K_{m+1}^\times)] + X^G$ in $[N_{K_m/F}
    (K_m^\times)]$. For $i<m$, then as before $Y_i^G$ is a
    complement of $[N_{K_{i+1}/F}(K_{i+1}^\times)] =
    (X+Y_{i+1}+\dots +Y_n)^{H_i}$ in $[N_{K_i/F}(K_i^\times)]$.
    Hence we have the statements establishing $\rank_{\Fp[G/H_i]}
    Y_i$.
\end{proof}

\section{Proofs of Theorem~3 and Corollary 3}\label{se:equivconds}

We shall first prove the first equality in
Theorem~\ref{th:equivconds} which says that
    \begin{equation*}
        m=i(K/F)=\min\{s \mid \xi_p\in N_{K/F}(K^\times)
        N_{K_{n-1}F}(K_s^\times)\}.
    \end{equation*}
In order to do so we shall calculate $(\root{p}\of{N_{K/F}
(\alpha)})^{\sigma-1}$, with a suitable $\alpha\in K^\times$,
in two ways. Then comparing our results we shall see that we
are indeed dealing with the equation
    \begin{equation*}
        E_i:\xi_p=N_{K/F}(\beta)N_{K_{n-1}/F}(\gamma),\quad\quad
        \beta\in K^\times,\ \gamma\in K^\times_i,\ 0\leq i<n
    \end{equation*}
and that our number $m=i(K/F)$ depends upon the smallest $i\in
\{-\infty,0,1,\dots,n-1\}$ such that $E_i$ is solvable for some
$\beta\in K^\times$ and $\gamma\in K_i^\times$. The following lemma
contains the key expression for the element $\left(\root{p}\of{N_
{K/F}(\alpha)}\right)^{\sigma-1}$.

\begin{lemma}\label{le:normcond}
    Suppose that $\alpha^{\sigma-1} = \gamma k^p$ with $\gamma
    \in K_i^\times$, $0\le i<n$, and $k\in K^\times$.  Suppose
    additionally that if $p=2$ then $n>1$.

    Then
    \begin{equation*}
        \left( \root{p}\of{N_{K/F}(\alpha)} \right) ^{\sigma-1} =
        N_{K/F}(k) \root{p}\of{N_{K/F}(\gamma)},
    \end{equation*}
    where
    \begin{equation*}
        \root{p}\of{N_{K/F}(\gamma)} = \left( N_{K_i/F}(\gamma)
        \right)^{p^{n-i-1}}.
    \end{equation*}
\end{lemma}

\begin{proof}
    First we claim that
    \begin{equation*}
        N_{K/F}(\alpha) = (k^p)^S\alpha^{p^n}\gamma^S,
    \end{equation*}
    where
    \begin{equation*}
        S := (p^{n}-1) + (p^{n}-2)\sigma + \dots +
        \sigma^{p^{n}-2} \in \Z[G].
    \end{equation*}
    Observe that
    \begin{align*}
        \alpha &= \alpha \\
        \alpha^\sigma &= k^p\alpha\gamma \\
        \alpha^{\sigma^2} &= \left( (k^p)^\sigma k^p \right)
        \alpha \left( \gamma \gamma^\sigma\right)
        \\
        \alpha^{\sigma^3} &= \left( (k^p)^{\sigma^2} (k^p)^\sigma
        k^p \right) \alpha \left( \gamma \gamma^{\sigma}
        \gamma^{\sigma^2} \right) \\
        &\dots \\
        \alpha^{\sigma^{p^n-1}} &= \left( \prod_{j=0}^{p^n-2}
        (k^p)^{\sigma^j} \right) \alpha \left( \prod_{j=0}^{p^n-2}
        \gamma^{\sigma^j}\right).
    \end{align*}
    Our result is then the product of the equations.

    Now $[\alpha]^{(\sigma-1)} = [\gamma]$, and because
    $[\gamma]\in [K_i^\times]$ and $[N_{K_i/F}(\beta)] =
    [\beta]^{(\sigma-1)^{p^i-1}}$ for $\beta\in K_i^\times$, we
    obtain $[\gamma]^{(\sigma-1)^{p^i}} = [1]$. Hence
    $[\alpha]^{(\sigma-1)^{p^i+1}} = [1]$.  Now $p^i+1 < p^n$
    unless $p=2$ and $n=1$, a case we have excluded. Hence
    $l(\alpha)<p^n$, whence $[N_{K/F}(\alpha)]=[1]$, and so
    $N_{K/F}(\alpha)\in K^{\times p}$. Therefore $\gamma^S \in
    K^{\times p}$ as well, and we may choose a $p$th root
    $\root{p}\of{\gamma^S} \in K^\times$. We then choose
    \begin{equation*}
        \root{p}\of{N_{K/F}(\alpha)} =
        k^S\alpha^{p^{n-1}}\root{p}\of{\gamma^S}.
    \end{equation*}

    (Because $(\root{p}\of{N_{K/F}(\alpha)})^{\sigma-1}$ does
    not depend upon the choice of a $p$th root of $N_{K/F}
    (\alpha)$ we see that we are free to make this choice.)

    Our next claim is that
    \begin{equation*}
        \left( \root{p}\of{\gamma^S} \right)^{\sigma-1} =
        \frac{(N_{K_i/F}(\gamma))^{p^{n-i-1}}}{\gamma^{p^{n-1}}}.
    \end{equation*}
    Let $L$ be the Galois closure of $K(\root{p}\of{\gamma})$ over
    $F$. Since $[\gamma]$ lies in the image of $\sigma-1$ on $J$,
    we have $[N_{K/F}(\gamma)]_F = [1]_F$.  Let $\hat\sigma$ be
    any pullback of $\sigma$ to $\Gal(L/F)$.  Then
    \begin{equation*}
        \root{p}\of{\gamma}^{(\hat\sigma^{p^n}-1)} =
        \root{p}\of{\gamma}^{(1+\hat\sigma+
        \dots+\hat\sigma^{p^n-1}) (\hat\sigma-1)} =
        \left(\root{p}\of{N_{K/F}(\gamma)} \right)^{
        (\hat\sigma-1)}=1.
    \end{equation*}
    (Observe that the equation is independent of the choice of the
    $p$th root of $N_{K/F}(\gamma)$.) Hence $\hat\sigma^{p^n}$
    leaves $\root{p}\of{\gamma}$ fixed. Now the field $L$ is
    generated over $K$ by all elements $\root{p}\of{\tilde\gamma}$,
    where $\tilde\gamma$ runs through all conjugate elements
    $\gamma^\tau$ for $\tau\in G=\Gal(K/F)$. Therefore
    $[N_{K/F}(\tilde \gamma)]_F=[1]_F$ for each such $\tilde\gamma$
    and the same argument as above shows that $\hat\sigma^{p^n}$
    leaves each $\root{p}\of{\tilde\gamma}$ fixed. Since
    $\hat\sigma$ restricted to $K$ is $\sigma$ we see that
    $\hat\sigma^{p^n}$ leaves every element of $K$ fixed. Hence
    $\hat\sigma^{p^n}$ leaves every element of $L$ fixed as well.
    Therefore $\hat \sigma^{p^n}=1\in \Gal(L/F)$.

    Set $\hat S = (p^n-1) + (p^n-2)\hat \sigma + \dots + \hat
    \sigma^{p^n-2}$, $\hat N = 1 + \hat \sigma + \dots + \hat
    \sigma^{p^n-1}$, and note that $\hat N = \hat N_1 \hat N_2$,
    where $\hat N_1 = 1 + \hat \sigma^{p^i} + \hat \sigma^{2p^i} +
    \dots + \hat \sigma^{(p^{n-i}-1)p^i}$ and $\hat N_2 = 1 + \hat
    \sigma + \dots + \hat \sigma^{p^i-1}$.  Further observe that
    $\hat N_1 \equiv N_{K/K_i}$ on $K^\times$, $\hat N_2 \equiv
    N_{K_i/F}$ on $K_i^\times$, and $(\hat\sigma-1)\hat S = \hat N -
    p^n$.

    We calculate $(\root{p}\of{\gamma^S})^{\sigma-1}$ in two cases.
    First assume that $\gamma\in K_0^\times$. Then $\gamma^S=\gamma^
    {p^n(p^{n}-1)/2}$ and since in the case $p=2$ we assume that
    $n\geq 2$ we see that $\gamma^S$ is a $p$th power of an element
    in $K_0^\times$ and therefore $(\root{p}\of{\gamma^S})^{\sigma-1}
    =1$ confirming our claim in this case. Next assume that $i>0$.
    Then we have
    \begin{align*}
        \left( \root{p}\of{\gamma^S} \right)^{\sigma-1} &=
        \left(\root{p}\of{\gamma}^{\hat S}\right)^{\sigma-1} =
        \frac{(\root{p}\of{\gamma})^{\hat
        N}}{(\root{p}\of{\gamma})^{p^n}} =
        \frac{\left(\root{p}\of{\gamma}^{\hat N_1}\right)^{\hat
        N_2}} {\gamma^{p^{n-1}}} \\ &=
        \frac{\left(\xi_p^c\gamma^{p^{n-i-1}}\right)^{\hat
        N_2}}{\gamma^{p^{n-1}}} =
        \frac{(N_{K_i/F}(\gamma))^{p^{n-i-1}}}{\gamma^{p^{n-1}}},
    \end{align*}
    where $\xi_p^c$ is a suitable $p$th root of $1$.

    Returning to $\root{p}\of{N_{K/F}(\alpha)}$, we may write
    \begin{align*}
        \left(\root{p}\of{N_{K/F}(\alpha)}\right)^{\sigma-1} &=
        k^{S(\sigma-1)} (\alpha^{p^{n-1}})^{\sigma-1}
        \left(\root{p}\of{\gamma^S}\right)^{\sigma-1} \\ &=
        k^{N-p^n}(\alpha^{\sigma-1})^{p^{n-1}}\left(\root{p}
        \of{\gamma^S}\right)^{\sigma-1} \\ &=\frac{N_{K/F}(k)}{k^
        {p^n}} \left( \gamma k^p\right)^{p^{n-1}}\frac{(N_
        {K_i/F}(\gamma))^{p^{n-i-1}}}{\gamma^{p^{n-1}}} \\ &=
        N_{K/F}(k)(N_{K_i/F}(\gamma))^{p^{n-i-1}}.
    \end{align*}
\end{proof}

\begin{proof}[of Theorem~{\rm\ref{th:equivconds}}]
    We have three equalities to establish, and we begin by showing
    that $m = \min \left\{s \ \vert \ \xi_p \in N_{K/F}(K^\times)
    N_{K_{n-1}/F}(K_{s}^\times)\right\} $.

    If $m=-\infty$ then $l(\delta)=1$ for $\delta$ an exceptional
    element.  Hence $[\delta]\in J^G$ and $[N_{K/F}(\delta)]_F\neq
    [1]_F$.  However, for all $f\in F^\times$, we have
    $[N_{K/F}(f)]_F = [1]_F$.  Hence $[\delta]\in J^G\setminus
    [F^\times]$.  Therefore, by the Fixed Submodule
    Lemma~\eqref{le:fixed}, $\xi_p \in N_{K/F}(K^\times)$. Going the
    other way, if $\xi_p \in N_{K/F}(K^\times)$, the Fixed Submodule
    Lemma~\eqref{le:fixed} tells us that there exists an exceptional
    element in $J^G$ and so $m=-\infty$. Hence the first equality of
    the Theorem holds when $m=-\infty$.

    Assume then that $m\ge 0$.  Consider $\alpha\in K^\times$ with
    $l(\alpha) < p^n$. By the Norm Lemma~\eqref{le:norm},
    $[N_{K/F}(\alpha)]_F \in \langle [a]_F\rangle$. It follows
    that $[N_{K/F}(\alpha)]_F \neq [1]_F$ if and only if
    $\root{p}\of{N_{K/F}(\alpha)}^{ \sigma-1}$ is a nontrivial
    $p$th root of unity, say $\xi^t$.

    Now assume that $[N_{K/F}(\alpha)]_F \neq [1]_F$ and
    $[\alpha]^{(\sigma-1)} = [\gamma]$, $\gamma\in K_i^\times$,
    $i<n$. Then by Lemma~\ref{le:normcond}, $\xi_p^t =
    N_{K/F}(k)(N_{K_i/F}(\gamma))^{p^{n-i-1}}$, and, by taking an
    appropriate power of $\xi_p^t$, we have
    \begin{equation*}
        \xi_p \in N_{K/F}(K^\times)(N_{K_i/F}
        (K_i^\times))^{p^{n-i-1}}.
    \end{equation*}
    Since there exists an exceptional element $\alpha$ with
    $[\alpha]^{(\sigma-1)} \in [K_m^\times]$, observe that
    $N_{K_{n-1}/F}(\gamma) = (N_{K_i/F}(\gamma))^{p^{n-i-1}}$ to
    conclude that $\xi_p \in N_{K/F}(K^\times)
    N_{K_{n-1}/F}(K_{m}^\times)$.  Hence the minimum $s$ is less
    than or equal to $m$.

    Going the other way, assume that $\xi_p = N_{K/F}(k)
    N_{K_s/F}(\gamma)^{p^{n-s-1}}$ for $k\in K$ and $\gamma\in
    K_s^\times$, $s<n$.  Then $1=N_{K/F}(k^p\gamma)$ and so by
    Hilbert~90 there exists $\delta\in K^\times$ with
    $\delta^{\sigma-1} = \gamma k^p$. Since $[\gamma]\in
    [K_s^\times]$ and $\sigma^{p^s}-1 \equiv (\sigma-1)^{p^s}$
    annihilates $[K_s^\times]$, we have $l(\delta) \le p^s+1 < p^n$.
    By Lemma~\ref{le:normcond}, $\root{p}
    \of{N_{K/F}(\delta)}^{(\sigma-1)} = N_{K/F}(k)(N_{K_s/F}
    (\gamma))^{p^{n-s-1}} = \xi_p$. Since $l(\delta) < p^n$ and
    $\root{p}\of{N_{K/F}(\delta)}^{(\sigma-1)}$ is a nontrivial
    $p$th root of unity, we use the equations above to deduce that
    $[N_{K/F}(\delta)]_F \neq [1]_F$.  Therefore by the definition
    of exceptionality, $m\le s$.

    We now establish the remaining two equalities.  For
    convenience, we set
    \begin{equation*}
        T := \left\{ t \ \ \vert \ \ \exists [\omega] \in J^{H_{t
        \dotplus 1}}, \ [N_{K/K_{t \dotplus 1}} (\omega)]_{K_{t
        \dotplus 1}} \neq [1]_{K_{t\dotplus 1}}\right\}
    \end{equation*}
    and
    \begin{equation*}
        S := \left\{ s \ \ \vert \ \ \xi_p \in N_{K/K_{s \dotplus
        1}} (K^\times) \right\}.
    \end{equation*}
    Observe that $n-1\in T$ because $\{0\} \neq X \subset J$ by
    Theorem~\ref{th:x} and $N_{K/K_n}(k)=k$ for each $k\in
    K^\times$, and $n-1\in S$ since $\xi_p \in F^\times\subset
    K^\times = N_{K/K_n}(K^\times)$. Hence the minima are
    well-defined.  It remains to show that $m = \min T = \min S$.

    To see that $m = \min T$, consider $t \in T$ with $t \leq n-2$
    such that there exists $[z] \in J^{H_{t\dotplus 1}}$ with
    $[N_{K/K_{t \dotplus 1}} (z)]_{K_{t \dotplus 1}} \neq [1]_{K_{t
    \dotplus 1}}$. By the Exact Sequence Lemma~\eqref{le:exact} we
    have $[N_{K/K_{t\dotplus 1}}(z)]_{K_{t\dotplus 1}} =
    [a_{t\dotplus 1}]_{K_{t\dotplus 1}}^s$ for some $s\in \Z$ with
    $s\not\equiv 0\bmod p$.  Then by the Norm Lemma~\eqref{le:norm}
    we see that $[N_{K/F}(z)]_F \neq [1]_F$. Now for $\delta$ an
    exceptional element of $K/F$, we have $l(\delta) =p^m + 1 \leq
    l(z)$, by Proposition \ref{pr:ldelta}, and hence $[z] \in
    J^{H_{t\dotplus 1}}$ implies that $l(\delta) = p^m + 1 \leq l(z)
    \leq p^{t\dotplus 1}$.  Hence $m \leq t$.  In the case $t =
    n-1$, again Proposition~\ref{pr:ldelta} gives $m<n$ and hence $m
    \leq t$. We conclude that $m \leq \min T$.

    For the other direction, observe that for $\delta$ an
    exceptional element of $K/F$, then $l(\delta)=p^m+1$ and
    therefore we have $[\delta] \in J^{H_{m\dotplus 1}}$.  By the
    Exact Sequence Lemma~\eqref{le:exact} the element
    $[N_{K/K_{m\dotplus 1}}(\delta)]_{K_{m\dotplus 1}}$ is contained
    in the subgroup of $K_{m\dotplus 1}^\times/K_{m\dotplus
    1}^{\times p}$ generated by $[a_{m\dotplus 1}]_{K_{m\dotplus
    1}}$.  By the Norm Lemma~\eqref{le:norm} and the remark
    following, we have $[N_{K/K_{m\dotplus
    1}}(\delta)]_{K_{m\dotplus 1}} \neq [1]_{K_{m\dotplus 1}}$. Hence
    $m \in T$ and $m \geq \min T$.  We conclude that $m = \min T$.

    Finally, we establish that $\min T = \min S$ by showing that
    $T=S$. Let $t \in T$ with $t \leq n-2$, and let $z$ satisfy $[z]
    \in J^{H_{t\dotplus 1}}$ and $[N_{K/K_{t \dotplus
    1}}(z)]_{K_{t\dotplus 1}} \neq [1]_{K_{t \dotplus 1}}$.  From
    the Fixed Submodule Lemma~\eqref{le:fixed}, part (2), we obtain
    $z^{(\sigma^{p^{t \dotplus 1}}-1)} = \lambda^p$ with
    $N_{K/K_{t\dotplus 1}} (\lambda) = \xi_p^\nu$ for some $\nu\in
    \Z$ with $\nu\not\equiv 0\bmod p$. Choosing an appropriate power
    of $\lambda$, we may assume that $\nu = 1$.  Hence $t \in S$.
    Since $n-1 \in T \cap S$, we have $T \subset S$.

    Conversely, suppose that $s \in S$ and $s\le n-2$ satisfies
    $\xi_p = N_{K/K_{s \dotplus 1}}(\lambda)$ for $\lambda \in
    K^\times$.  We have $1 = N_{K/K_{s \dotplus 1}}(\lambda^p)$, and
    by Hilbert 90 we see that there exists $\delta \in K^\times$
    such that $\delta^{\sigma^{p^{s \dotplus 1}}-1} = \lambda^p$.
    Hence $[\delta] \in J^{H_{s\dotplus 1}}$, and again using the
    Fixed Submodule Lemma~\eqref{le:fixed} and its proof we see that
    $[N_{K/K_{s\dotplus 1}} (\delta)]_{K_{s \dotplus 1}} \neq
    [1]_{K_{s\dotplus 1}}$.  Hence $s \in T$. Since $n-1 \in T \cap
    S$, we have $S \subset T$.
\end{proof}

\begin{proof}[of Corollary~{\rm\ref{co:equivconds}}]
    Assume that $\xi_p \in F$ and, if $p=2$, then either $n>1$ or
    $-1 \in N_{K/F}(K^\times)$. Consider the case $i(K/F)\ge 0$ first,
    and choose $j\in \{0, \dots, i(K/F)\}$. Then
    \begin{align*}
        i(K/K_j) &= \min \{ t \mid \xi_p \in N_{K/
        K_{j+t+1}}(K^\times) \} \\ &= \min \{ s \mid \xi_p \in
        N_{K/K_{s+1}}(K^\times) \} - j \\ &= i(K/F) - j.
    \end{align*}
    These equalities are clear, except possibly in the case when
    $p=2$, $i(K/F) = n-1$, and $j = n-1$.  In this case we must also
    verify that $i(K/K_{n-1})$ is well defined and is equal to $0$.

    Now if $n=1$ then Proposition~\ref{pr:excexist} tells us that
    $i(K/K_{n-1})$ is defined, since we have assumed that $-1 \in
    N_{K/F}(K^\times)$. In order to show that $i(K/K_{n-1})$ is well
    defined when $n>1$, it is sufficient to show that the set
    \begin{equation*}
        \{ \delta \in K^\times \mid [N_{K/K_{n-1}}
        (\delta)]_{K_{n-1}} \neq [1]_{K_{n-1}} \}
    \end{equation*}
    is not empty.

    Suppose the set is empty: $[N_{K/K_{n-1}}(K^\times)]_{K_{n-1}}=
    [1]_{K_{n-1}}$. Using elements $a_i \in K^\times_i$, $0 \le i <
    n$, as in Proposition~\ref{pr:subgen} we have $K =
    K_{n-1}(\sqrt{a_{n-1}})$ and $-a_{n-1} \in
    N_{K/K_{n-1}}(K^\times)$. Thus $[-a_{n-1}]_{K_{n-1}} =
    [1]_{K_{n-1}}$ and we have $a_{n-1} = -\nu^2$ for some $\nu \in
    K^\times_{n-1}$. Therefore, using the hypothesis $n>1$, we see
    that $[N_{K_{n-1}/F}(a_{n-1})]_F = [1]_F$, contrary to our
    choice of $a_{n-1}$.

    Because
    \begin{equation*}
        [\sqrt{a_{n-1}}]^{(\sigma-1)} = [\sqrt{a_{n-1}}]^{
        (\sigma+1)} = N[\sqrt{a_{n-1}}]\in [K_{n-1}^\times],
    \end{equation*}
    we see that $i(K/K_{n-1})\le 0$.  If $i(K/K_{n-1}) = -\infty$
    then by the Exact Sequence Lemma~\eqref{le:exact} we have $-1\in
    N_{K/K_{n-1}}(K^\times)$ and therefore $i(K/F)\le n-2$, a
    contradiction of our hypothesis.  Hence $i(K/K_{n-1})=0$.

    Now assume that $i(K/F) = -\infty$. Then $\xi_p \in
    N_{K/F}(K^\times)$ and by \cite[Theorem~3]{A} we see that $\xi_p
    \in N_{K/K_j}(K^\times)$ as well.  Therefore $i(K/K_j)=-\infty$. If
    $0 \le i(K/F) < j$, then similarly $\xi_p \in
    N_{K/K_j}(K^\times)$, and we again conclude that $i(K/K_j) =
    -\infty$.

    Finally, for $j<n$ the cyclic extension $K_j/F$ embeds
    into the cyclic extension $K_{j+1}/F$.  By \cite[Theorem~3]{A},
    we have $\xi_p\in N_{K_j/F}(K_j^\times)$, and hence
    $i(K_j/F)=-\infty$.
\end{proof}

\section*{Acknowledgements}

Andrew Schultz would like to thank Ravi Vakil for his encouragement
and direction in this and all other projects. All of the authors are
very grateful to the anonymous referee for reading this paper very
carefully and providing a number of valuable suggestions concerning
the exposition.

\end{document}